\documentclass[11pt,a4paper,twoside,bibtotoc,final]{scrartcl}
\usepackage{a4wide}
\usepackage{amsfonts}
\usepackage{amsmath}
\usepackage{amssymb}
\usepackage{amsthm}
\usepackage{mathtools}
\usepackage[numbers]{natbib} 
\usepackage{booktabs}
\usepackage[title]{appendix}
\usepackage[Algorithm]{algorithm}
\usepackage{algorithmic}
\usepackage{pgfplots}
\pgfplotsset{compat=newest}
\usepackage{graphicx}
\usepackage{soul}
\usepackage{color}
\usepackage{colortbl}
\usepackage{multirow}
\usepackage{boxedminipage}
\usepackage{enumerate}
\usepackage[caption=false,font=footnotesize]{subfig}
\usepackage{rotating}
\usepackage{hyperref}
\usepackage{todonotes}

\newlength\fheight
\newlength\fwidth
\newcommand{\N}{\ensuremath{\mathbb{N}}}
\newcommand{\E}{\ensuremath{\mathbb{E}}}
\newcommand{\T}{\ensuremath{\mathbb{T}}}

\newcommand{\Z}{\ensuremath{\mathbb{Z}}}

\newcommand{\R}{\ensuremath{\mathbb{R}}}

\newcommand{\C}{\ensuremath{\mathbb{C}}}
\newcommand{\I}{\ensuremath{I}}
\newcommand{\J}{\ensuremath{\mathrm{J}}}
\renewcommand{\P}{\ensuremath{\mathcal{P}}}
\newcommand{\X}{\mathcal{X}}
\newcommand{\bigO}{\ensuremath{\mathcal{O}}}

\newcommand{\ii}{\textnormal{i}}
\newcommand{\e}{\textnormal{e}}

\newcommand{\zb}[1]{\ensuremath{\boldsymbol{#1}}}

\renewcommand{\ln}{\mathrm{ln\,}}

\newcommand{\boldk}{{\ensuremath{\boldsymbol{k}}}}
\newcommand{\boldi}{{\ensuremath{\boldsymbol{i}}}}

\newcommand{\boldl}{{\ensuremath{\boldsymbol{l}}}}
\newcommand{\boldh}{{\ensuremath{\boldsymbol{h}}}}

\newcommand{\boldx}{{\ensuremath{\boldsymbol{x}}}}
\newcommand{\boldy}{{\ensuremath{\boldsymbol{y}}}}
\newcommand{\boldz}{{\ensuremath{\boldsymbol{z}}}}

\newcommand{\boldI}{{\ensuremath{\boldsymbol{I}}}}

\newcommand{\boldzero}{{\ensuremath{\boldsymbol{0}}}}
\newcommand{\boldone}{{\ensuremath{\boldsymbol{1}}}}
\newcommand{\setT}{{\ensuremath{\mathcal{T}}}}

\newcommand{\err}{\ensuremath{\operatorname{err}}}

\DeclareMathOperator*{\argmax}{arg\,max}

\renewcommand{\d}{\,\mathrm{d}}

\definecolor{light-gray}{gray}{0.75} 
\sethlcolor{light-gray}

\newtheorem{theorem}{Theorem}[section]
\newtheorem{lemma}[theorem]{Lemma}
\newtheorem{remark}[theorem]{Remark}
\newtheorem{generalisation}[theorem]{Generalisation}
\newtheorem{definition}[theorem]{Definition}
\newtheorem{example}[theorem]{Example}
\newtheorem{corollary}[theorem]{Corollary}
\newtheorem{proposition}[theorem]{Proposition}

\makeatletter
\def\imod#1{\allowbreak\mkern10mu({\operator@font mod}\,\,#1)}
\makeatother

\numberwithin{equation}{section}
\numberwithin{table}{section}
\numberwithin{figure}{section}

\newcommand{\bend}{\hspace*{0ex} \hfill \hbox{\vrule height
    1.5ex\vbox{\hrule width 1.4ex \vskip 1.4ex\hrule  width 1.4ex}\vrule
    height 1.5ex}}

\long\def\symbolfootnote[#1]#2{\begingroup%
\def\thefootnote{\fnsymbol{footnote}}\footnote[#1]{#2}\endgroup}
\newcommand{\sspan}{\textnormal{span}}

\newcommand{\OO}[1]{\mathcal{O}\left(#1\right)}

\renewcommand{\mathbf}[1]{\ensuremath{\boldsymbol{#1}}}
\renewcommand{\textbf}[1]{{\ensuremath{\boldsymbol{#1}}}}
\renewcommand{\thefootnote}{\fnsymbol{footnote}}
\pgfplotscreateplotcyclelist{MR1LOF}{%
dashed, mark=diamond, mark options={solid, scale=1.5}\\%
dashed, mark=o, mark options={solid, scale=1.5}\\%
densely dotted, mark=triangle, mark options={solid, scale=1.5}\\%
densely dotted, mark=square, mark options={solid, scale=1.25}\\%
dashdotted,mark=x, mark options={solid, scale=1.5}\\%
dashdotted, mark=+, mark options={solid, scale=1.5}\\%
loosely dashed, mark=square, mark options={solid,scale=0.75}\\%
loosely dashdotted, mark=asterisk, mark options={solid,scale=0.5}, opacity=0.5\\%
dasdotdotted, every mark/.append style={solid},mark=star\\%
densely dashdotted,every mark/.append style={solid, fill=gray},mark=diamond*\\%
}%

\usetikzlibrary{patterns}
\makeatletter %
\pgfdeclarepatternformonly[\LineSpace,\tikz@pattern@color]{my north east lines}{\pgfqpoint{-1pt}{-1pt}}{\pgfqpoint{\LineSpace}{\LineSpace}}{\pgfqpoint{\LineSpace}{\LineSpace}}%
{
    \pgfsetcolor{\tikz@pattern@color} %
    \pgfsetlinewidth{0.4pt}
    \pgfpathmoveto{\pgfqpoint{0pt}{0pt}}
    \pgfpathlineto{\pgfqpoint{\LineSpace + 0.1pt}{\LineSpace + 0.1pt}}
    \pgfusepath{stroke}
}
\makeatother %
\newdimen\LineSpace
\tikzset{
    line space/.code={\LineSpace=#1},
    line space=3pt
}
\title{The uniform sparse FFT with application to PDEs with random coefficients}

\date{September 1, 2022}
\author{
Lutz K\"ammerer\footnotemark[1], Daniel Potts\footnotemark[2], Fabian Taubert\footnotemark[3]}

\hypersetup{pdfauthor={Lutz K\"ammerer,Daniel Potts,Fabian Taubert},
          pdftitle={The uniform sparse FFT with application to PDEs with random coefficients},
          pdfsubject={},
          pdfcreator = {pdflatex},
          plainpages=false,
          pdfstartview=FitH,       
          pdfview=FitH,            
          pdfpagemode=UseOutlines, 
          bookmarksnumbered=true, 
          bookmarksopen=false,     
          bookmarksopenlevel=0,   
          colorlinks=true,       
          linkcolor=black,
          citecolor=black,
          urlcolor=black}

\newif\ifshowextendedpaperversion
\showextendedpaperversionfalse

\begin{document}

\maketitle

\begin{abstract}
\small
We develop the uniform sparse Fast Fourier Transform (usFFT), an efficient, non-intrusive, adaptive algorithm for the solution of elliptic partial differential equations with random coefficients. The algorithm is an adaption of the sparse Fast Fourier Transform (sFFT), a dimension-incremental algorithm, which tries to detect the most important frequencies in a given search domain and therefore adaptively generates a suitable Fourier basis corresponding to the approximately largest Fourier coefficients of the function. The usFFT does this w.r.t.\ the stochastic domain of the PDE simultaneously for multiple fixed spatial nodes, e.g., nodes of a finite element mesh. The key idea of joining the detected frequency sets in each dimension increment results in a Fourier approximation space, which fits uniformly for all these spatial nodes. This strategy allows for a faster and more efficient computation due to a significantly smaller amount of samples needed, than just using other algorithms, e.g., the sFFT for each spatial node separately. We test the usFFT for different examples using periodic, affine and lognormal random coefficients in the PDE problems.

\medskip

{\small%

}
\end{abstract}\

\footnotetext[1]{
  Chemnitz University of Technology, Faculty of Mathematics, 09107 Chemnitz, Germany\\
  kaemmerer@mathematik.tu-chemnitz.de
}
\footnotetext[2]{
  Chemnitz University of Technology, Faculty of Mathematics, 09107 Chemnitz, Germany\\
  potts@mathematik.tu-chemnitz.de
}
\footnotetext[3]{
  Chemnitz University of Technology, Faculty of Mathematics, 09107 Chemnitz, Germany\\
  fabian.taubert@mathematik.tu-chemnitz.de
}
\medskip

\ifshowextendedpaperversion
\tableofcontents
\newpage
\fi

\section{Introduction}\label{sec:intro}
Parametric operator equations have gained significant attention in recent years. In particular, partial differential equations with random coefficients play an important role in the study of uncertainty quantification, e.g., \cite{BoRaSchw17,KaKaKu20,KaKuSlo20}. Therefore, the numerical solution of these equations and how to compute them in an efficient and reliable way has become more and more important.

In this work, we consider the parametric, elliptic problem of finding $u: D_{\boldx} \times D_{\boldy} \rightarrow \R$ such that for every $\boldy \in D_{\boldy}$ there holds
\begin{equation}
\begin{split}
-\nabla \cdot \left( a(\boldx,\boldy) \nabla u(\boldx,\boldy) \right) &= f(\boldx) \qquad\qquad\boldx \in D_{\boldx},\, \boldy \in  D_{\boldy},\\
u(\boldx,\boldy) &= 0 \quad\qquad\qquad\boldx \in \partial D_{\boldx},\, \boldy \in D_{\boldy},
\end{split}\label{eq:pde} 
\end{equation}
describing the diffusion characteristics of inhomogeneous materials and therefore being called diffusion equations with the random diffusion coefficients $a$. Here, $\boldx = (x_j)_{j=1}^{d_{\boldx}} \in D_{\boldx}$ is the spatial variable in a bounded Lipschitz domain $D_{\boldx} \subset \R^{d_{\boldx}}$, typically with spatial dimension $d_{\boldx}=1,2$ or $3$, and $\boldy = (y_j)_{j=1}^{d_{\boldy}} \in D_{\boldy}$ is a high-dimensional random variable with $D_{\boldy} \subset \R^{d_{\boldy}}$. For the remainder of this paper, we use $d$ instead of $d_{\boldy}$ to simplify notations. The differential operator $\nabla$ is always used w.r.t.\ the spatial variable $\boldx$ and the one-dimensional random variables $y_j$ are assumed to be i.i.d. with a prescribed distribution. 

A common way to define the random coefficient $a$ is via
\begin{align}\label{eq:random_coefficient_1}
a(\boldx,\boldy) = a_0(\boldx) + \sum_{j=1}^{d} \Theta_j(\boldy) \,\psi_j(\boldx),
\end{align}
where $a_0$ and $\psi_j$ are assumed to be uniformly bounded on $D_{\boldx}$. This model is commonly used with the stochastic domain $D_{\boldy} = [\alpha,\beta]^{d}$, typically with $[\alpha,\beta] = [-1,1]$ or $[-\frac12, \frac12]$. The $\Theta_j(\boldy)$ can also be interpreted as random variables itself and are usually chosen to have expectation value $\E[\Theta_j(\boldy)]=0$, such that $\E[a(\boldx,\cdot)] = a_0(\boldx)$ holds and the terms of the sum model the stochastic fluctuations.

Often, the model \eqref{eq:random_coefficient_1} is in an affine fashion, using $\Theta_j(\boldy) = y_j$ for all $j=1,...,d$. This so-called affine model is considered in many works on parametric differential equations with random coefficients, e.g., \cite{CoDeSchw10b, KuSchwSl12, Schw13, ZhaHuHoLiYa14, DiLGSchw14, KuSchwSl15, BaCoDa16,  DiKuLG16, BaCoGi17, GaHeSchw18, NgNu21uniform}.
The so-called periodic model using $\Theta_j(\boldy) = \frac{1}{\sqrt{6}}\sin(2\pi y_j)$ has been recently studied in \cite{KaKuSlo20, KaKaKu20}. For $y_j$ uniformly distributed on $[-\frac12,\frac12]$ each, these $\Theta_j$ are then distributed according to the arcsine distribution on $[-1,1]$. It turned out that this model is also worth to be considered in addition to the affine model. Further, this model yields some advantages for our new approach due to its periodicity w.r.t.\ the random variables, as we will see later.

The second type of the random coefficient $a$, that is also used in many recent works, e.g., \cite{GrKuNi13, CheHoYaZha13, BaCoDVGi17, BaCoDu17, BaCoMi18, NgNu21lognormal}, is the so-called lognormal form
\begin{align*}
&a(\boldx,\boldy) = a_0(\boldx) + \exp(b(\boldx,\boldy)), & b(\boldx,\boldy) = b_0(\boldx) + \sum_{j=1}^{d} y_j \,\psi_j(\boldx).
\end{align*}
Here, the random variables $y_j$ are typically normally distributed, i.e., $y_j \sim \mathcal{N}(0,1)$, and hence $D_{\boldy} = \R^{d}$. The numerical analysis as well as the computation of approximations for this model is more difficult, but also arises more often from real applications. A more detailed overview on parametric and stochastic PDEs can be found, e.g., in \cite[Sec.~1]{CoDe15}.

In this paper, we design a numerical method for solving the aforementioned problems. To be more precise, we will compute approximations of the solutions $u(\boldx,\boldy)$ using trigonometric polynomials. A Fourier approach on ordinary differential equations, i.e., $d_x=1$, with high-dimensional random coefficients has already been presented in \cite{BoKaePo20}. There, a dimension-incremental method, the so-called \textit{sparse Fast Fourier Transform (sFFT)}, cf.\ \cite{PoVo14}, was used to detect the most important frequencies $\boldk$ and corresponding approximations of the Fourier coefficients $c_{\boldk}(u)$ of the solution $u(\boldx,\boldy)$. These values can be used to compute an approximation of the solution $u$ or other quantities of interest as, e.g., the expectation value $\E[u]$. Further, the frequencies and Fourier coefficients can be used to gain detailed information about the influence of the random variables $y_j$ on the solution $u$ and their interaction with each other. 

In this work, we present a non-intrusive approach based on the main idea of the algorithm developed in \cite{BoKaePo20}. The main difference is, that we do not include the spatial variable $\boldx$ in the Fourier approach and therefore only apply the sFFT w.r.t.\ the random variable $\boldy$. Therefore, the sFFT only needs samples of the function values of $u$ for fixed $\boldy$, which can be computed by using any suitable, already available differential equation solver. In consequence, we are not restricted to particular spatial domains $D_{\boldx}$ or spatial dimensions $d_x$. To be more precise, we consider a finite set $\setT_G \subset D_{\boldx}$ with finite cardinality $|\setT_G| = G < \infty,$ as spatial discretization and aim for approximations of the functions $$u_{\boldx_g}(\boldy) \coloneqq u(\boldx_g,\boldy)$$ for each $\boldx_g \in \setT_G$. Also note that we might need to apply a suitable periodization w.r.t.\ $\boldy$ first if the function $u_{\boldx_g}$ is not already periodic.

Unfortunately, we would need to apply such a pointwise approximation algorithm, like in our case sFFT, then $G$ times separately, resulting in an unnecessary huge increase in the number of samples used and therefore, since each sample implies a call of the underlying, probably expensive differential equation solver, also in computation time of the algorithm. Hence, we develop a modification of the sFFT to overcome this problem and compute the approximations of the functions $u_{\boldx_g}$ within one call of the new algorithm. In particular, our so-called \textit{uniform sparse Fast Fourier Transform (usFFT)} combines some candidate sets between each dimension-incremental step, which allows to use the same sampling nodes $\boldy$ for each point $\boldx_g \in \setT_G$ in the next step. This strategy manages to keep the number of used samples in a reasonable size and hence decreases the computation time drastically compared to $G$ applications of the sFFT algorithm itself. We summarize this in the following Theorem:

\begin{theorem}\label{thm:complexities}
Let the sparsity parameter $s\in\N$, a frequency candidate set $\Gamma\subset\Z^d$, $|\Gamma|<\infty$, the amount $G\in\N$ and a failure probability $\delta\in(0,1)$ be given. Moreover, we define
$N_{\Gamma} \coloneqq \max_{j=1,...,d} \lbrace \max_{\boldk\in\Gamma} k_j - \min_{\boldl \in \Gamma} l_j \rbrace$. 
Then, there exists a randomized sampling strategy based on the random rank-1 lattice approach in \cite{KaKrVo20} generating a set $S$ of sampling locations with cardinality 
\begin{align}\label{eq:sample_compl}
|S|\in\bigO \left( d\, s\, \max(s, N_{\Gamma})\, \log^2 \frac{d\, s\, G\, N_{\Gamma}}{\delta} + \max(sG, N_\Gamma)\, \log \frac{d\,s\,G}{\delta} \right)
\end{align}
such that the following holds.

Consider $G$ arbitrary multivariate trigonometric polynomials
$p^{(g)}(\boldy):=\sum_{\boldk\in \I_g}\hat{p}_{\boldk}^{(g)}\e^{2\pi\ii\boldk\cdot\boldy}$, $g=1,\ldots,G$, where we assume $\I_g\subset\Gamma$, $|\I_g|\le s$ and $\min_{\boldk\in \I_g}|\hat{p}_{\boldk}^{(g)}|>0$ for each $g=1,\ldots,G$.
We generate a random set $S$ via this sampling strategy. Then, with probability at least $1-\delta$ it holds that
\begin{itemize}
\setlength{\itemindent}{2em}
\item all frequencies $\boldk\in \I_g$ as well as
\item all Fourier coefficients $\hat{p}_{\boldk}^{(g)}$, $\boldk\in \I_g$,
\end{itemize}
of all multivariate trigonometric polynomials $p^{(g)}$, $g=1,\ldots,G$, can be reconstructed from their values at the sampling locations in $S$.

The simultaneous identification of all the frequencies and the computation of all the Fourier coefficients can be realized by a combination of Algorithm~\ref{algo:sfft_general} and a modification of the approach presented in \cite{KaKrVo20} in the role of Algorithm $\operatorname{A}$. The suggested method has a computational complexity of
\begin{align*}
\bigO \left( d^2\, s^2 \,G^2 \,N_{\Gamma}\, \log^3 \frac{d \,s\, G\, N_{\Gamma}}{\delta}\right)
\end{align*}
with probability at least $1-\delta$ as well as $\bigO \left( d^2 \,s^3 \,G^2 \,N_{\Gamma} \,\log^3 \frac{d \,s \,G \,N_{\Gamma}}{\delta}\right)$ in the worst case.
\end{theorem}

\begin{remark}\label{rem:complex}
Note that \eqref{eq:sample_compl} in Theorem \ref{thm:complexities} does not state anything about the sampling complexity of Algorithm \ref{algo:sfft_general}, but the amount of sampling locations $\vert S \vert$. We need samples of all trigonometric polynomials $p^{(g)}(\boldy), g=1,\ldots,G,$ at these sampling nodes $\boldy \in S$, so the necessary amount of samples in the classical sense is $G$ times larger. However, we aim for an Algorithm, where the $G$ samples $p^{(g)}(\boldy)$ for a fixed $\boldy \in S$ for all $g=1,\ldots,G$ are obtained by just a single call of some (probably expensive) black box sampling method. In all of our numerical examples in Section \ref{sec:numerics}, the computation time for this sampling procedure tremendously outweighs the pure computation time of the remaining steps of Algorithm \ref{algo:sfft_general}, which we referred to as computational complexity in Theorem \ref{thm:complexities}. Hence, we stress on the fact, that the computational complexity is not the main focus of this complexity result, but the amount of sampling nodes.
\end{remark}

The proof of the Theorem is given in Appendix \ref{app:proof}. While Theorem \ref{thm:complexities} is stated for trigonometric polynomials $p^{(g)}$ only, the algorithm can be used on the above mentioned periodic or periodized functions $u_{\boldx_g}(\boldy)$ to compute the support and values of the approximately largest Fourier coefficients of the functions with some thresholding technique, realized by the parameters $s_\mathrm{local}$ and $\theta$ in Algorithm \ref{algo:sfft_general}, as well, which is also the key idea when applying the sFFT for function approximation in \cite{PoVo14, KaPoVo17, KaKrVo20}. More generally spoken, we could even consider $G$ different periodic functionals $F_g(\boldy)$ and approximate them with the same approach we are about to present here. Moreover, Theorem \ref{thm:complexities} does not assume the frequency sets $\I_g$ to share any frequencies $\boldk$, i.e., these sets could even be pairwise disjoint in the worst case scenario. Obviously, this will not be the case in our examples later on as the functions $u_{\boldx_g}(\boldy)$ and $u_{\boldx_{\tilde{g}}}(\boldy)$ are probably very similar for $\boldx_g$ and $\boldx_{\tilde{g}}$ close to each other due to the smoothness of the solution $u(\boldx,\boldy)$. Hence, the given complexities, especially the quadratic dependency on $G$ of the computational complexity, are very pessimistic and should really be seen as a worst case estimate.

The crucial advantage of the presented approach is the efficient and adaptive choice of the frequency set performed by the underlying sFFT. Most of the approaches in the aforementioned works are based on certain (tensorized) basis functions \cite{CoDeSchw10b, CheHoYaZha13, BaCoDa16, BaCoGi17, BaCoDVGi17, BaCoDu17, BaCoMi18, KaKaKu20}, Quasi-Monte Carlo methods \cite{KuSchwSl12, Schw13, GrKuNi13, CheHoYaZha13, DiLGSchw14, KuSchwSl15, DiKuLG16, GaHeSchw18, KaKuSlo20, NgNu21uniform, NgNu21lognormal} or collocation methods \cite{CheHoYaZha13, ZhaHuHoLiYa14} and often assume the particular involved basis functions, weights, index sets or kernels needed to be known, chosen or computed in advance. A common example are compressed sensing techniques, see, e.g., \cite[Ch.~7]{AdBruWeb22} or \cite{FoRa13} and the references therein for an overview, where the considered index sets need to be chosen in advance. The adaptivity of our algorithm circumvents this step and therefore provides much more freedom in finding a good sparse approximation. Also note that our method aims for the approximation of the solution $u(\boldx,\boldy)$ directly instead of, e.g., just a high-dimensional quadrature. As an example and for additional information, we refer to \cite{KuNu18} as a short and general introduction to Quasi-Monte Carlo methods, which is also one of the most common approaches.

Another suitable approach is to use an efficient deterministic sampling strategy which guarantees the reconstruction of each $s$-sparse signal supported on the $d$-dimensional candidate set $\Gamma$, so-called deterministic sparse Fourier Transform. Such a method allows to use the same samples for approximating all signals $u_{\boldx_g}$, $g=1,\ldots, G$.
Several works like \cite{Iw13} already provide applicable multivariate sparse Fourier transform results, but the resulting complexities $\bigO (d^4s^2)$ (up to logarithmic factors) in a deterministic setting and $\bigO (d^4s)$ for a random variant scale suboptimal in the dimension $d$. Another fully deterministic method is stated in \cite{Mo17}. There, the construction of the sampling sets suffers from some minor restrictions on the considered frequency set, which result in a slightly better scaling sampling complexity $\bigO (d^3s^2N)$ and computational complexity $\bigO (d^3s^2N^2)$ (both again up to logarithmic factors) with $N$ the side length of the considered cube in frequency domain. Finally, in \cite{GrIwKaVo20} multivariate sparse Fourier transforms are presented, where the corresponding complexities scale again quadratically in $s$ for a deterministic version and linearly for a Monte Carlo version, while the dimension $d$ only enters linearly. Unfortunately, this is only true if the considered candidate sets do not scale exponentially w.r.t.\ the dimension $d$. Otherwise, the size $M$ of the reconstructing rank-1 lattices used also scales exponentially in $d$ and logarithmic factors like $\log^{11/2} M$ in the complexities then result again in a suboptimal dimension scaling.

Our usFFT is highly adaptive, since it only needs an arbitrary candidate set $\Gamma$ and selects the important frequencies $\boldk$ in this search domain on its own. The cardinality $|\Gamma|$ of this candidate set is not as crucial as for other approaches like mentioned above, since the number of used samples and the computation time suffer only mildly from larger candidate sets. Again, we stress on the fact that we may also extract additional information about the influence and the interactions of the random variables $\boldy_j$ from the output of the usFFT. For instance, we detect a maximum of only $4$ simultaneously active dimensions in the detected frequencies in our numerical examples, i.e., the detected frequency vectors $\boldk$ have at most $4$ non-zero components with $d = 10$ or even $d =20$. Another main advantage of our algorithm is the non-intrusive and parallelizable behavior. As already mentioned, the usFFT uses existing numerical solvers of the considered differential equation. We can use suitable, reliable and efficient solvers with no need to re-implement them. Further, the different samples needed in each sampling step can be computed on multiple instances. This parallelization allows to reduce the computation time even further and makes a higher number of used samples less time consuming.

The remainder of the paper is organized as follows: 

In Section \ref{sec:prerequ} we set up some notation and assumptions and briefly explain the key idea of the sFFT algorithm. Section \ref{sec:usfft} is devoted to the explanation of the usFFT as well as some periodizations required for the affine and lognormal cases. Finally, in Section \ref{sec:numerics} we test the new algorithm on different examples using periodic, affine, and lognormal random coefficients and investigate the computed approximations under different aspects.

The MATLAB\textsuperscript{\textregistered} source code of the algorithm as well as demos for our numerical examples can be downloaded from \url{https://mytuc.org/fyfw}.

\section{Prerequisites}
\label{sec:prerequ}

We consider the PDE problem \eqref{eq:pde}. Note that we always assume $f$ to be independent of the random variable $\boldy$ and zero boundary conditions just for simplicity and to preserve clarity. Our algorithm (up to some minor changes) may also be applied for right-hand sides $f(\boldx,\boldy)$ as well as non-zero Dirichlet boundary conditions $u(\boldx,\boldy) = h(\boldx,\boldy)$ for all $\boldx \in \partial D_{\boldx}$.

\subsection{Problem setting}
\label{subsec:pde}

The weak formulation of our problem reads: Given $f \in H^{-1}(D_{\boldx})$, for every $\boldy \in D_{\boldy}$, find $u(\cdot,\boldy) \in H_0^1(D_{\boldx})$, such that
\begin{align*}
\int_{D_{\boldx}} a(\boldx,\boldy) \nabla u(\boldx,\boldy) \cdot \nabla v(\boldx) \d\boldx = \int_{D_{\boldx}} f(\boldx)v(\boldx)\d\boldx \quad \forall v \in H_0^1(D_{\boldx}).
\end{align*}
As usual, $H_0^1(D_{\boldx})$ denotes the subspace of the $L_2$-Sobolev space $H^1(D_{\boldx})$ with vanishing trace on $\partial D_{\boldx}$ and $H^{-1}(D_{\boldx})$ denotes the dual space of $H_0^1(D_{\boldx})$.
We say, that the diffusion coefficient $a: D_{\boldx} \times D_{\boldy} \rightarrow \R$ fulfills the uniform ellipticity assumption, if there exist two constants $a_{\min} \in \R$ and $a_{\max} \in \R$, such that
\begin{align}\label{eq:UEA}
&0 < a_{\min} \leq a(\boldx,\boldy) \leq a_{\max} < \infty & \forall \boldx \in D_{\boldx},\, \forall \boldy \in D_{\boldy}.
\end{align}
Then, the Lax-Milgram Lemma ensures, that the problem \eqref{eq:pde} possesses a unique solution $u(\cdot,\boldy) \in H_0^1(D_{\boldx})$ for every fixed $\boldy \in D_{\boldy}$, satisfying the a priori estimate
\begin{align*}
\sup_{\boldy \in D_{\boldy}} \| u(\cdot,\boldy) \|_{H_0^1(D_{\boldx})} \leq \frac{1}{a_{\min}} \| f \|_{H^{-1}(D_{\boldx})}.
\end{align*}
Some further basic information and results on approximation and smoothness of the solution $u$ of high-dimensional parametric PDEs can be found in \cite[Sec.~1 and 2]{CoDe15}. Additionally, we also refer to the general results on best $n$-term approximations given in \cite[Sec.~3.1]{CoDe15}, since our Fourier approach fits in this particular framework as well.

In order to compute an approximation of the solution $u_{\boldx_g} \coloneqq u(\boldx_g,\cdot)$ at a given point $\boldx_g \in D_{\boldx}$ using the dimension-incremental method explained below, we need samples of $u_{\boldx_g}$ for a lot of sampling nodes $\boldy$. We aim for a non-intrusive approach and therefore use a finite element method to solve the problem \eqref{eq:pde} for a given $\boldy \in D_{\boldy}$. A similar approach is used e.g.\ in \cite{NgNu21uniform, NgNu21lognormal}, where the finite element method is used to solve the PDE for any $\boldy_{\mathfrak{u}}$ with $\mathfrak{u} \subset \N$ and $(\boldy_{\mathfrak{u}})_j = y_j$ for $j\in \mathfrak{u}$ and $0$ otherwise. The corresponding approximations of the so-called $\mathfrak{u}$-truncated solution are then used for their particular method aswell. In our case, we just evaluate the finite element solution $\check{u}(\cdot,\boldy)$ at the given point $\boldx_g \in D_{\boldx}$. In particular, instead of the finite element method, any differential equation solver would fit, that is capable of computing the value $u(\boldx_g,\boldy)$ for given $\boldx_g$ and $\boldy$. Hence, we also refer to this sampling method as black box sampling later on.

Note that we will use the finite element solution $\check{u}$ also as an approximation of the true solution $u$, when we test the accuracy of our computed approximation $u^{\texttt{usFFT}}$ in Section \ref{sec:numerics}. In detail, we have
\begin{align*}
\text{err}(u,u^{\texttt{usFFT}}) \leq \text{err}(u,\check{u}) + \text{err}(\check{u},u^{\texttt{usFFT}}),
\end{align*}
where $\text{err}(\cdot,\cdot)$ is a suitable metric, symbolizing the error. So while we only investigate the second term $\text{err}(\check{u},u^{\texttt{usFFT}})$ in our numerical tests later, the first term includes other error sources as the modeling, e.g., by a dimension truncation of infinite-dimensional random coefficient $a$, or the error coming from the finite element approximation itself. For a particular example of this, we refer to the detailed error analysis for the periodic model mentioned in Section~\ref{sec:intro}, that is given in \cite[Sec.~4]{KaKaKu20}.

\subsection{The dimension-incremental method for $s$-sparse periodic functions}
\label{subsec:method}

The following dimension-incremental method was presented in \cite{PoVo14}. The aim of this algorithm is to determine the non-zero Fourier coefficients $\hat{p}_{\boldk} \in \C, \, \boldk \in \I,$ of a multivariate trigonometric polynomial
\begin{align*}
p(\boldy) = \sum_{\boldk \in \I} \hat{p}_{\boldk} \exp(2\pi\ii\boldk\cdot\boldy)
\end{align*}
with unknown frequency set $\I \subset \Z^d,\, |\I| < \infty,$ based on samples of the polynomial $p$.
Obviously, $p$ is a periodic signal and its domain is the $d$-dimensional torus $\T^d$, $\T\simeq [0,1)$.

The goal is not only to calculate the nonzero Fourier coefficients $\hat{p}_{\boldk}$ but also, and more important, to detect the frequencies $\zb k$ out of a possibly huge search domain $\Gamma\subset\Z^d$ belonging to the nonzero Fourier coefficients. In particular, we define the set $$\operatorname{supp}\hat{p}:=\{\boldk\in \Gamma\colon \hat{p}_\boldk\neq 0\} $$ and call the cardinality $|\operatorname{supp}\hat{p}|$ the sparsity of $p$.

First, we introduce some further notation. We consider a given search domain $\Gamma \subset \Z^d, \, |\Gamma| < \infty,$ that should be large enough to contain the unknown frequency set ${\I \subset \Gamma}$. We denote the projection of a frequency $\boldk \coloneqq (k_1,...,k_d) \in \Z^d$ to the components $\boldi \coloneqq (i_1,...,i_m) \in \lbrace \iota \in \lbrace 1,...,d \rbrace^m : \iota_t \not = \iota_{t^{\prime}} \text{ for } t \not = t^{\prime} \rbrace$ by $\P_{\boldi}(\boldk) \coloneqq (k_{i_1},...,k_{i_m}) \in \Z^m$. Correspondingly, we define the projection of a frequency set $\I \subset \Z^d$ to the components $\boldi$ by $\P_{\boldi}(\I) \coloneqq \lbrace (k_{i_1},...,k_{i_m}): \boldk \in \I \rbrace$. Using these notations, the general approach is the following:

\subsubsection*{Sketch of dimension-incremental reconstruction}

\begin{enumerate}
\item Compute the first components of the unknown frequency set from sampling values, i.e., determine a set $\I^{(1)} \subset \P_1(\Gamma)$, such that $\P_1(\text{supp } \hat{p}) \subset \I^{(1)}$ holds.

\item For dimension increment step $t=2,...,d$, i.e., for each additional dimension:
\begin{enumerate}
\item Compute the $t$-th components of the unknown frequency set from sampling values, i.e., determine a set $\I^{(t)} \subset \P_t(\Gamma)$, such that $\P_t(\text{supp } \hat{p}) \subset \I^{(t)}$ holds.

\item Construct a suitable sampling set $\X^{(1,...,t)} \subset \mathbb{T}^d,\, |\X^{(1,...,t)}| \ll |\Gamma|$, which allows to detect those frequencies from the set $(\I^{(1,...,t-1)} \times \I^{(t)}) \cap \P_{(1,...,t)}(\Gamma)$ belonging to non-zero Fourier coefficients $\hat{p}_{\boldk}$.

\item Sample the trigonometric polynomial $p$ along the nodes of the sampling set $\X^{(1,...,t)}$.

\item Compute the Fourier coefficients $\tilde{\hat{p}}_{(1,...,t),\boldk},\, \boldk \in (\I^{(1,...,t-1)} \times \I^{(t)}) \cap \P_{(1,...,t)}(\Gamma)$.

\item Determine the non-zero Fourier coefficients from $\tilde{\hat{p}}_{(1,...,t),\boldk},\, \boldk \in (\I^{(1,...,t-1)} \times \I^{(t)}) \cap \P_{(1,...,t)}(\Gamma)$ and obtain the set $\I^{(1,...,t)}$ of detected frequencies. The $\I^{(1,...,t)}$ index set should be equal to the projection $\P_{(1,...,t)} (\text{supp }\hat{p})$.
\end{enumerate}

\item Use the set $\I^{(1,...,d)}$ and the computed Fourier coefficients $\tilde{\hat{p}}_{(1,...,d),\boldk},\, \boldk \in \I^{(1,...,d)}$ as an approximation for the support $\text{supp } \hat{p}$ and the Fourier coefficients $\hat{p}_{\boldk},\, \boldk \in \text{supp } \hat{p}$.
\end{enumerate}

Note that this method can also be used for the numerical determination of the  approximately largest Fourier coefficients 
$$c_{\boldk}(f):=\int_{\T^d}f(\boldy)\e^{-2\pi\ii\boldk\cdot\boldy}\mathrm{d}\boldy,\quad \boldk \in \I,$$ of suffciently smooth periodic signals $f$ using a suitable thresholding technique, cf.
\cite{PoVo14,KaPoVo17,KaKrVo20}.

The proposed approach includes the construction of suitable sampling sets in step~2b.
To this end, one assumes that an upper bound $s\ge|\operatorname{supp}\hat{p}|$ is known and one constructs the sampling sets $\mathcal{X}^{(1,\ldots,t)}$ such that the Fourier coefficients $\tilde{\hat{p}}_{(1,...,t),\boldk}$ computed in step~2d
are randomly projected ones. Due to that projection one may observe cancellations with the effect that one misses active frequencies. For that reason, one repeats the computation of the projected Fourier coefficients for a number $r$ of random projections and then one takes the union.

Of course, there exist different methods for the computation of the projected Fourier coefficients. The algorithm works with any sampling method, which computes Fourier coefficients on a given frequency set. Preferable sampling sets combine the four properties:
\begin{itemize}
\setlength{\itemsep}{0cm}
\setlength{\parskip}{0cm}
\item relatively low number of sampling nodes (sampling complexity),
\item stability,
\item efficient construction methods for the sampling set,
\item fast Fourier transform like algorithms in order to compute the projected Fourier coefficients.
\end{itemize}
Especially due to the last point, we will call the dimension-incremental method the \textit{sparse Fast Fourier Transform (sFFT)} from now on. A quick sketch of the sampling techniques used in \cite{PoVo14,KaPoVo17,KaKrVo20} as well as the sample complexity and computational complexity of the sFFT using these methods are given in Appendix \ref{app:r1l}.

\section{The uniform sparse FFT}\label{sec:usfft}

Up to now, the sFFT algorithm is a suitable tool in order to compute an approximation of the solution $$u_{\boldx_g}(\boldy) \coloneqq u(\boldx_g,\boldy)=\sum_{\boldk\in\Z^d}c_\boldk(u_{\boldx_g})\e^{2\pi\ii\boldk\cdot\boldy}\approx\sum_{\boldk\in \I_{\boldx_g}}c_\boldk^{\texttt{sFFT}}(u_{\boldx_g})\e^{2\pi\ii\boldk\cdot\boldy}$$ for a single $\boldx_g$. When considering a whole set of points $\boldx_g \in \setT_G$, $\lvert \setT_G \rvert = G$, we have to call the existing method $G$ times. But multiple, independent calls of the sFFT result in different, adaptively determined sampling sets $\X^{(1,...,t)}$ in step~2b. Hence, we cannot guarantee that the solutions of the differential equation from one run of the algorithm can be utilized in another one. So we really need $G$ full calls of the sFFT including all sampling computations and therefore end up with unnecessary many samples, even when using the sample efficient rank-1 lattice (R1L) approaches. Remember, that sampling means solving the differential equation with a call of the underlying differential equation solver, that might be very expensive in computation time. Therefore, we now modify the dimension-incremental method, such that we can work on the set $\setT_G$ and one call of the algorithm computes approximations of the most important Fourier coefficients $c_{\boldk}(u_{\boldx_g})$, $\boldk\in \I_{\boldx_g}$, for each $g=1,...,G$, including a clever choice of the sampling nodes $\boldy$.

\subsection{Expanding the sFFT}\label{subsec:usFFT}

The full method is stated in \autoref{algo:sfft_general}. We force the dimension-incremental method to select a set $\I\subset\Gamma\subset\Z^d$ containing the frequencies of the $s$ approximately largest Fourier coefficients $c_{\boldk}(u_{\boldx_g})$ for each $\boldx_g \in \setT_G$. To this end, we compute the set of detected frequencies $\I_{\boldx_g}^{(1,...,t)}$ for each $\boldx_g$ in each dimension-increment $t$, but afterwards we form the union of these sets $\bigcup_{g=1}^{G} \I_{\boldx_g}^{(1,...,t)}$, which will be the set of detected frequencies $\I^{(1,...,t)}$, that is given to the next dimension-incremental step $t+1$. 

Now, we start each iteration with a larger frequency candidate set $(\I^{(1,...,t-1)} \times \I^{(t)}) \cap \P_{(1,...,t)}(\Gamma)$, which is suitable for all $\boldx_g \in \setT_G$.
This way, the first $t$ components of the elements of the sampling set $\X^{(1,...,t)}$ are the same for each $\boldx_g$ and the random part, which causes the specific random projection of the Fourier coefficients, can be chosen equally for each $\boldx_g$ without disturbing the algorithm.
Therefore, we can now take advantage of the fact, that our underlying differential equation solver can evaluate the solutions $u(\boldx,\boldy)$ for a given $\boldy$ for multiple values of $\boldx$ in the domain $D_{\boldx}$. Accordingly, we only need to solve the differential equation once for each sampling node $\boldy$ and still get all the sampling values $u_{\boldx_g}(\boldy)$ for all $\boldx_g$ in our finite set $\setT_G$. Note that this also holds for the one-dimensional detections in steps~1 and~2a. Also, we might need to interpolate or approximate, if some of the values $u_{\boldx_g}(\boldy)$, $\boldx_g\in \setT_G$ and fixed $ \boldy \in \X^{(1,...,t)}$, are not directly given by the differential equation solver. Obviously, the larger candidate sets $(\I^{(1,...,t-1)} \times \I^{(t)}) \cap \P_{(1,...,t)}(\Gamma)$, resulting from the union of the sets $\I_{\boldx_g}^{(1,...,t-1)}$, $g=1,\ldots,G$, and the union of the sets $\I_{\boldx_g}^{(t)}$, $g=1,\ldots,G$, will also result in larger sampling sets. The overall increase of sampling locations considered is still very reasonable, cf.\ Theorem \ref{thm:complexities}. The computational complexity suffers a bit harder from these modifications, but is not as important as the amount of sampling locations, cf.\ Remark \ref{rem:complex}. We could also think of further thresholding methods to cut the number of frequencies back to the sparsity parameter $s$ at the end of each dimension-incremental step or at least at the end of the whole algorithm. In this work, we will not do this, but take a look at the total number of frequencies in the output of the algorithm in relation to the sparsity parameter $s$, cf.\ Remark \ref{rem:growth}.

Overall, the proposed method is now capable of computing approximations
\begin{align*}
u_{\boldx_g}^{\texttt{usFFT}}(\boldy) \coloneqq \sum_{\boldk\in \I}c_\boldk^{\texttt{usFFT}}(u_{\boldx_g})\e^{2\pi\ii\boldk\cdot\boldy}
\end{align*}
for all nodes $\boldx_g, g=1,...,G.$ Please note that step~2 does not provide $(c_\boldk^{\texttt{usFFT}}(u_{\boldx_g}))_{\boldk\in \I}$ but only approximations of $(c_\boldk^{\texttt{usFFT}}(u_{\boldx_g}))_{\boldk\in \tilde{J}_{d,1,g}}$, where $\tilde{J}_{d,1,g}\subsetneq \I^{(1,\ldots,d)}=\I$ holds in general.
In order to compute all the Fourier coefficients $c_\boldk^{\texttt{usFFT}}(u_{\boldx_g})$, $\boldk\in \I$ and $g=1,\ldots, G$, we propose an additional approximation in step~3. For this approximation, the user can apply a suitable approach of his choice. The used method should just compute an approximation of the projection to the already determined space of trigonometric polynomials $\sspan\{\exp(2\pi\ii\boldk\cdot\circ)\colon \boldk\in I\}$, cf., e.g., \cite{KuoMiNoNu19, Kae17, KaUlVo19} for different possible sampling approaches.

We will call this modified version of the sFFT the \textit{uniform sFFT} or short \textit{usFFT} from now on, where \textit{uniform} is meant w.r.t.\ the discrete set of points $\setT_G$. The main difference to the sFFT algorithms from \cite{PoVo14,KaPoVo17,KaKrVo20} are the loops over $g$ and the corresponding unions of the frequency sets. Possible choices for Algorithm~$\operatorname{A}$ and the approximation approach used in step~3 of Algorithm~\ref{algo:sfft_general} are a random rank-1 lattice approach and a multiple rank-1 lattice approach, respectively, which actually leads to Theorem~\ref{thm:complexities}, cf.\ Appendix \ref{app:r1l} and Appendix \ref{app:proof}.

\begin{algorithm}[p]
	\caption{The usFFT on a set $\setT_G$}\label{algo:sfft_general}
  \begin{small}
	\begin{tabular}{p{1.2cm}p{2.1cm}p{11cm}}
		Input:  %
		& $\Gamma\subset\Z^d$ \hfill & search space in frequency domain, candidate set for $\I$ \\
		& $u(\cdot,\cdot)$ & PDE solution $u$ as black box (function handle) \\
        & $\setT_G$ & discrete set containing the points $\boldx_g,\, g=1,...,G$\\
		& $s,s_\mathrm{local}\in\N$ & sparsity parameters, $s\leq s_\mathrm{local}$ \\  
		& Algorithm $\operatorname{A}$ & {\bfseries efficient algorithm} $\operatorname{A}$ that guarantees the identification of the frequency support of each $s_\mathrm{local}$-sparse trigonometric polynomial with high probability, cf.\ Section \ref{subsec:method}, and computes the Fourier coefficients \\
		& $\theta\in\R^+$ & absolute threshold \\
		& $r\in\N$ & number of detection iterations%
	\end{tabular}
	\vspace{-0.6em}
	\begin{algorithmic}
		\item[(Step 1 \& 2a)] [Single frequency component identification]
		\STATE {\bfseries for} $t:=1,\ldots,d$ {\bfseries do}
		\STATE \hspace{1em} Set $K_t:=\max(\mathcal{P}_t(\Gamma))-\min(\mathcal{P}_t(\Gamma))+1$.
		\STATE \hspace{1em} Set $\I^{(t)}:=\emptyset$.
		\STATE \hspace{1em} {\bfseries for} $i:=1,\ldots,r$ {\bfseries do}
		\STATE \hspace{2em} Choose $y_j'\in\mathbb{T}$, $j\in\{1,\ldots,d\}\setminus\{t\}$ uniformly at random.
		\STATE \hspace{2em} Set $\boldy^{(\ell)}:=\left(y_1^{(\ell)},\ldots, y_d^{(\ell)}\right)^\top$, $y_j^{(\ell)}:=
		\begin{cases}
		\ell/K_t,& j=t,\\
		y_j',& j\neq t,
		\end{cases}$ \quad
		for all $\ell=0,\ldots,K_t-1$.
		\STATE \hspace{2em} {\bfseries for} $g:=1,\ldots,G$ {\bfseries do}
		
   		\STATE \hspace{3em} Compute $\tilde{\hat{u}}_{t,k_t,g}:=\frac{1}{K_t} \sum_{\ell=0}^{K_t-1} u\left(\boldx_g,\boldy^{(\ell)}\right) \,\mathrm{e}^{-2\pi\mathrm{i} \ell k_t / K_t}$, $k_t\in\mathcal{P}_t(\Gamma)$, via FFT.
  		\STATE \hspace{3em} Set $\I^{(t)}:=\I^{(t)}\cup\{k_t\in\mathcal{P}_t(\Gamma)\colon \tilde{\hat{u}}_{t,k_t,g} \text{ is among the largest } s_\mathrm{local} \text{ (in absolute value)}$
		\STATE \hspace{3em} \qquad\qquad\qquad\qquad $\text{ elements of } \{\tilde{\hat{u}}_{t,j,g}\}_{j\in\mathcal{P}_t(\Gamma)} \text{ and } \vert\tilde{\hat{u}}_{t,k_t,g}\vert \geq \theta
		\}$.   
		\STATE \hspace{2em} {\bfseries end for} $g$
		\STATE \hspace{1em} {\bfseries end for} $i$
		\STATE {\bfseries end for} $t$
		\item[(Step 2)] [Coupling frequency components identification]
		\STATE {\bfseries for} $t:=2,\ldots,d$ {\bfseries do}
		\STATE \hspace{1em} If $t<d$, set $\tilde r:=r$ and $\tilde s:=s_\mathrm{local}$, otherwise $\tilde r:=1$ and $\tilde s:=s$. 
		\STATE \hspace{1em} Set $\I^{(1,\ldots,t)}:=\emptyset$.
		\STATE \hspace{1em} \parbox[t]{\linewidth-\algorithmicindent-\algorithmicindent-2em}{Generate a sampling set $\mathcal{X}\subset\mathbb{T}^t$ for $ \J_t:=(\I^{(1,\ldots,t-1)} \times \I^{(t)})\cap\mathcal{P}_{(1,\ldots,t)}(\Gamma)$ that allows for the application of  Algorithm $\operatorname{A}$.}
		\STATE \hspace{1em} {\bfseries for} $i:=1,\ldots,\tilde r$ {\bfseries do}
		\STATE \hspace{2em} Choose components $y_{t+1}',\ldots,y_d'\in\mathbb{T}$ of sampling nodes uniformly at random.
		\item[(Step 2b)]
		\STATE \hspace{2em} Set $\mathcal{X}_{t,i}:=\{\boldy:=(\boldsymbol{\tilde{y}},y_{t+1}',\ldots,y_d')\colon \boldsymbol{\tilde{y}}\in\mathcal{X} \}\subset\mathbb{T}^d$.
		\item[(Step 2c)]
		\STATE \hspace{2em} \parbox[t]{\linewidth-\algorithmicindent-\algorithmicindent-2em}{Sample $u$ along the nodes of the sampling set $\mathcal{X}_{t,i}$ for every $\boldx_g$.}
		\STATE \hspace{2em} {\bfseries for} $g:=1,\ldots,G$ {\bfseries do}
		\item[(Step 2d)]
		\STATE \hspace{3em} \parbox[t]{\linewidth-\algorithmicindent-\algorithmicindent-3em}{Apply Algorithm~$\operatorname{A}$ to obtain the support $\tilde{\J}_{t,i,g}\subset \J_t$, $|\tilde{\J}_{t,i,g}|\le \tilde{s}$, of frequencies belonging to the at most $\tilde{s}$ largest Fourier coefficients, each larger than $\theta$ in absolute value, using the sampling values~$u(\boldx_g,\boldy_j)$, $\boldy_j\in\mathcal{X}_{t,i}$.}	
		\item[(Step 2e)]
		\STATE \hspace{3em} Set $\I^{(1,\ldots,t)}:=\I^{(1,\ldots,t)}\cup \tilde{\J}_{t,i,g}$.
     	\STATE \hspace{2em} {\bfseries end for} $g$
		\STATE \hspace{1em} {\bfseries end for} $i$
		\STATE {\bfseries end for} $t$
   \end{algorithmic}
   \end{small}
\end{algorithm}

\begin{algorithm*}
  \begin{small}
  \ContinuedFloat
  \caption{continued.}
  \begin{algorithmic}
		\item[(Step 3)] [Computation of Fourier coefficients]
		\STATE Generate a suitable sampling set $\mathcal{Y}\subset\T^d$.
		\STATE Sample $u$ along the nodes of the sampling set $\mathcal{Y}$ for every $\boldx_g$. 
		\STATE {\bfseries for} $g:=1,\ldots,G$ {\bfseries do}
		\STATE \hspace{1em} Compute the corresponding Fourier coefficients $\left(\tilde{\hat{u}}_{(1,\ldots,d),\boldk,g}\right)_{\boldk\in \I^{(1,\ldots,d)}}$ 
		\newline \phantom{STATE} by the means of the samples $\left(u(x_g,\boldy) \right)_{\boldy\in\mathcal{Y}}$ and a suitable algorithm.
		\STATE {\bfseries end for} $g$
		\STATE Set $\tilde{\I}:=\I^{(1,\ldots,d)}$
	\end{algorithmic}

	\begin{tabular}{p{1.3cm}p{2.0cm}p{11.0cm}}
		Output: & $\tilde{\I}\subset\Gamma\subset\Z^d$ & set of detected frequencies\\
		& $\boldsymbol{\tilde{\hat{u}}}_{g}\in\C^{\vert \tilde{\I}\vert}$ & corresponding Fourier coefficients for all $\boldx_g$, where each $|\tilde{\hat{u}}_{(1,\ldots,d),\boldk,g}| \geq \theta$ for at least one $\boldx_g$ \\
	\end{tabular}
  \end{small}
\end{algorithm*}

\subsection{Periodization}

The usFFT allows us to reconstruct a frequency set $\I$ and approximations $c_{\boldk}^{\texttt{usFFT}}(u_{\boldx_g})$ of the corresponding Fourier coefficients $c_{\boldk}(u_{\boldx_g})$ for each $\boldx_g \in \setT_G$. Unfortunately, this approach requires the function $u(\boldx,\boldy)$ to be $1$-periodic w.r.t.\ $\boldy$ in each stochastic dimension $d$. 

Since the right-hand side $f(\boldx)$ does not depend on $\boldy$ in our considerations, the random coefficient $a$ is the only given function involving the random variable $\boldy$ in the problem \eqref{eq:pde}. In periodic models, we use the random coefficient \eqref{eq:random_coefficient_1} with 1-periodic functions $\Theta_j(\boldy)$. Hence, the random coefficient $a(\boldx,\boldy)$ is 1-periodic and thus the solution $u(\boldx,\boldy)$ is also 1-periodic w.r.t.\ each component of $\boldy$. Therefore, we can apply the usFFT directly for this model without any further considerations. 

In order to apply the usFFT when using the affine and lognormal models, we need to apply a suitable periodization first, since the random coefficient $a$ and therefore the solution $u$ are not periodic in general. Note that we assume the random variable to be uniformly distributed in the affine case, i.e.,\ $\boldy \sim \mathcal{U}([\alpha,\beta]^{d})$, and standard normally distributed in the lognormal case, i.e.,\ $\boldy \sim \mathcal{N}(0,1)^{d}$ and recall $\T \simeq [0,1)$.

\subsubsection{Affine case}
\label{subsubsec:affine}
We consider the in $\tilde{\boldy}$ 1-periodic function
\begin{align*}
&\tilde{u}: D_{\boldx} \times \T^{d} \longrightarrow \R \\
&\tilde{u}(\boldx,\tilde{\boldy}) \coloneqq u(\boldx,\varphi(\tilde{\boldy})),
\end{align*}
with $\varphi$ being some suitable transformation function, i.e.,
\begin{align*}
\varphi: \T^{d} \longrightarrow D_{\boldy} = [\alpha,\beta]^{d}.
\end{align*}
With this approach, the usFFT is able to compute approximations of the functions \linebreak ${\tilde{u}_{\boldx_g} \coloneqq \tilde{u}(\boldx_g,\cdot)}$ for each $\boldx_g \in \setT_G$. We want $\varphi$ to act component-wise on the random variable, i.e., $\varphi(\tilde{\boldy}) \coloneqq \left( \varphi_j(\tilde{y}_j) \right)_{j=1}^{d}$. Further, we assume, that these mappings $\varphi_j$ fulfill the assumptions
\begin{enumerate}
\item[(A1)] Each $\varphi_j$ is continuous, i.e., $\varphi_j \in C(\T)$ for each $j=1,...,d$.
\item[(A2)] It holds $\varphi_j(0) = \varphi_j(1) = \alpha$ and $\varphi_j(\frac12) = \beta$ for each $j=1,...,d$.
\item[(A3)] Each $\varphi_j$ is symmetric, i.e., $\varphi_j(\frac12-\tilde{y}) = \varphi_j(\frac12+\tilde{y})$ for $\tilde{y} \in [0, \frac12]$ and for each $j=1,...,d$.
\item[(A4)] Each $\varphi_j$ is strictly monotonously increasing in $[0,\frac12]$ for each $j=1,...,d$.
\end{enumerate}
With these restrictions we ensure, that $\varphi$ is bijective w.r.t.\ the interval $[0,\frac12]^{d}$. Hence, we define the inverse mapping $\varphi^{-1}(\boldy): [\alpha,\beta]^{d} \rightarrow [0,\frac12]^{d}$.

With this inverse mapping, we are now able to compute approximations of the functions $u_{\boldx_g}(\boldy)$ via
\begin{align}\label{eq:eval_formula}
u_{\boldx_g}^{\texttt{usFFT}}(\boldy) \coloneqq \tilde{u}_{\boldx_g}^{\texttt{usFFT}}(\varphi^{-1}(\boldy)) = \sum_{\boldk \in \I} c_{\boldk}^{\texttt{usFFT}}(\tilde{u}_{\boldx_g})\,\e^{2\pi\ii\boldk\cdot \varphi^{-1}(\boldy)},
\end{align}
with the finite index set $\I$ and the approximated Fourier coefficients $c_{\boldk}^{\texttt{usFFT}}(\tilde{u}_{\boldx_g})$ from the usFFT applied to the functions $\tilde{u}_{\boldx_g}$, $g\in\mathcal{T}_G$.

In this work, we always consider the tent transformation, cf.\ \cite{LiHi03,SuNuCo14,CoKuNuSu16}, for each $\varphi_j$, i.e.,
\begin{subequations}\label{eq:tent}
\begin{align}
&\varphi_j: \T \longrightarrow [\alpha,\beta],& \varphi_j(\tilde{y}) = \beta - \left|(\beta-\alpha)\left(1 - 2\tilde{y}\right)\right|,\\
&\varphi_j^{-1}: [\alpha,\beta] \longrightarrow \left[0,\frac12\right],& \varphi_j^{-1}(y) = \frac{y-\alpha}{2(\beta-\alpha)}.
\end{align}
\end{subequations}
Although this transformation mapping fulfills the assumptions (A1) - (A4), it might not be the most favorable choice in specific applications due to its lack of smoothness. Smoother periodizations, e.g., \cite[Sec.~2.2.2]{BoKaePo20}, might yield better approximation results in specific situations due to the faster decay of the Fourier coefficients of $\tilde{u}$. On the other hand, the linear structure  of the tent transformation on the interval $[0,\frac12]$ allows some simplifications later on.

\subsubsection{Lognormal case}
As in the affine case, we need a suitable, periodic transformation mapping ${\varphi: \T^{d} \rightarrow D_{\boldy} = \R^{d}}$ to receive a periodization $\tilde{u}(\boldx,\tilde{\boldy})$. Again, we choose the same functions in each stochastic dimension, so the same $\varphi_j$ for all $j=1,...,d$, but this time $\varphi_j$ will consist of two separate steps. First, we consider the transformation
\begin{align*}
&\tau_1: \left(-\frac12,\frac12\right) \longrightarrow \R,& \tau_1(\breve{y}) \coloneqq \sqrt{2} \, \text{erf}^{-1}(2\breve{y}),\\
&\tau_1^{-1}: \R \longrightarrow \left(-\frac12,\frac12\right),& \tau_1^{-1}(y)= \frac12 \, \text{erf}\left(\frac{y}{\sqrt{2}}\right),
\end{align*}
with the error function
\begin{align*}
\text{erf}(y) \coloneqq \frac1{\sqrt{\pi}} \int_{-y}^{y} \e^{-t^2} \d t, \quad x \in \R.
\end{align*}
For further information on this transformation, see \cite{NaPo18}. This mapping $\tau_1$ seems like the ideal choice when talking about random variables $\boldy \sim \mathcal{N}(0,1)$, since the error function $\text{erf}(y)$ is closely related to its cumulative distribution function $\Phi$. In detail, it holds
\begin{align*}
\Phi(y) = \frac{1}{2} \left( 1+ \text{erf} \left( \frac{y}{\sqrt{2}} \right)\right).
\end{align*}
The so-called inversion method in stochastic simulation describes, that the cumulative distribution function $\Phi$ and its inverse $\Phi^{-1}$ map random variables, distributed according to $\Phi$, to uniformly distributed random variables on $[0,1]$ and the other way around, cf.\ \cite[Sec.~II.2]{Devroye86}. Thus, our transformation $\tau_1$ maps uniformly distributed random variables $\breve{y} \sim \mathcal{U}(-\frac12,\frac12)$ to normally distributed random variables $y \sim \mathcal{N}(0,1)$ and is therefore a great generalization when moving forward from uniformly distributed random variables.

The second part is a suitable periodization $\tau_2: \T \rightarrow (-\frac12,\frac12)$. We choose a similar approach as in the affine case and use a shifted tent transformation
\begin{align*}
&\tau_{2,\Delta}: \T \longrightarrow \left[-\frac12,\frac12\right], & \tau_{2,\Delta}(\tilde{y}) = \begin{cases}
-\frac12 - 2(\tilde{y}-\Delta) & 0 \leq \tilde{y} < \Delta \\
-\frac12 + 2(\tilde{y}-\Delta) & \Delta \leq \tilde{y} < \frac12+\Delta \\
+\frac32 - 2(\tilde{y}-\Delta) & \frac12+\Delta \leq \tilde{y} < 1
\end{cases} \\
&\tau_{2,\Delta}^{-1}: \left[-\frac12,\frac12\right] \longrightarrow \left[\Delta,\frac12+\Delta\right], & \tau_{2,\Delta}^{-1}(\breve{y}) = \frac{\breve{y}}{2} + \Delta + \frac14
\end{align*}
with shift $\Delta > 0$. We need this shift, since we cannot apply the transformation $\tau_1$ if we have $\tau_{2,\Delta}(\tilde{y})=\pm \frac12$ due to the poles there. Shifting with a suitable $\Delta$ ensures, that the deterministic part of the sampling set $\X$ does not contain components equal to $\Delta$ or $\frac12 + \Delta$. The randomly chosen values from the interval $[0,1]$ for the other components will not be equal to $\Delta$ or $\frac12 + \Delta$ almost surely too. Hence, the sampling set $\X$ in \autoref{algo:sfft_general} does not contain any nodes with any component equal to $\Delta$ or $\frac12+\Delta$ almost surely. 

\begin{figure}[tb]
	\centering	
	\subfloat[transformation mapping $\tau_1$]{
    		\setlength\fwidth{3.5cm}
\definecolor{mycolor1}{rgb}{0.00000,0.44700,0.74100}%
\begin{tikzpicture}
\begin{axis}[%
width=0.951\fwidth,
height=\fheight,
at={(0\fwidth,0\fheight)},
scale only axis,
xmin=-0.5,
xmax=0.5,
xtick={-0.5, -0.25, 0, 0.25, 0.5},
ymin=-4,
ymax=4,
axis background/.style={fill=white},
xmajorgrids,
ymajorgrids,
legend style={legend cell align=left, align=left, draw=white!15!black}
]
\addplot [color=mycolor1]
  table[row sep=crcr]{%
-0.499998999999	-4.75342410672217\\
-0.499712999713	-3.4436111621475\\
-0.499090999090999	-3.11845778805354\\
-0.498146998146998	-2.90215901492863\\
-0.496880996880997	-2.73500112755026\\
-0.495291995291995	-2.59656862601528\\
-0.493377993377993	-2.47713948338564\\
-0.491134991134991	-2.37120867523186\\
-0.488559988559989	-2.27543207824312\\
-0.485648985648986	-2.18755505575401\\
-0.482397982397982	-2.10600365105718\\
-0.478802978802979	-2.02963118659652\\
-0.474858974858975	-1.95755671819763\\
-0.47056097056097	-1.88910277362205\\
-0.465902965902966	-1.82372228389616\\
-0.460878960878961	-1.76097828611582\\
-0.455482955482956	-1.70051504432572\\
-0.44970794970795	-1.64202848683308\\
-0.443544943544944	-1.58524749657006\\
-0.436986936986937	-1.5299620369527\\
-0.43002393002393	-1.47596927755299\\
-0.422646922646923	-1.42310349724412\\
-0.414844914844915	-1.37120785477632\\
-0.406605906605907	-1.32014027903688\\
-0.397915897915898	-1.26976541333699\\
-0.388760888760888	-1.21996479087924\\
-0.379124879124879	-1.17062326455088\\
-0.368989868989869	-1.12162889177372\\
-0.358334858334858	-1.07286813845955\\
-0.347137847137847	-1.0242349686743\\
-0.335371835371835	-0.975612900384792\\
-0.323006823006823	-0.926884792188683\\
-0.31000881000881	-0.877928761011903\\
-0.296337796337796	-0.828611272441827\\
-0.281944781944782	-0.778778107403135\\
-0.266773766773767	-0.728263229852571\\
-0.250753750753751	-0.676863603771307\\
-0.233798733798734	-0.624342722384585\\
-0.215796715796716	-0.570399795262213\\
-0.196603696603697	-0.514657171159033\\
-0.176018676018676	-0.456594340304251\\
-0.153755653755653	-0.395479981547761\\
-0.129371629371629	-0.330189551463414\\
-0.102107602107602	-0.258806171658515\\
-0.07041757041757	-0.177437365941077\\
-0.0297625297625297	-0.0746729372635659\\
0.0817825817825817	0.206455865294844\\
0.113158613158613	0.287561036818209\\
0.14000964000964	0.358484560689494\\
0.163978663978664	0.423346226374264\\
0.185839685839686	0.484091929929543\\
0.206037706037706	0.541846016662792\\
0.224858724858725	0.5973367842124\\
0.242499742499742	0.651071217962681\\
0.25910575910576	0.703428931562124\\
0.274785774785775	0.754700926594508\\
0.289624789624789	0.805120025689076\\
0.303692803692804	0.854885763716948\\
0.317046817046817	0.90416792378783\\
0.32973282973283	0.95310999884761\\
0.341791841791842	1.00184943672109\\
0.353256853256854	1.05050434683974\\
0.364158864158864	1.09919671412975\\
0.374522874522874	1.14803467167148\\
0.384371884371884	1.19712912038144\\
0.393725893725894	1.2465890670181\\
0.402602902602903	1.29652638820347\\
0.411018911018911	1.34705606101972\\
0.418988918988919	1.39830278464927\\
0.426526926526926	1.45040308635926\\
0.433643933643934	1.5034923324138\\
0.440350940350941	1.55772640097906\\
0.446657946657947	1.61327807963344\\
0.452573952573952	1.67034009370522\\
0.458107958107958	1.72913976468645\\
0.463266963266963	1.78992435165019\\
0.468058968058968	1.85300205658393\\
0.472490972490973	1.91873361514198\\
0.476568976568976	1.98754488669693\\
0.48029898029898	2.05996340499445\\
0.483686983686984	2.13665536856252\\
0.486737986737987	2.2184520354713\\
0.489457989457989	2.30647651198644\\
0.491850991850992	2.40217260380455\\
0.493922993922994	2.50764087888795\\
0.495677995677996	2.62582168819222\\
0.497122997122997	2.76147990417486\\
0.498262998262998	2.92235050555215\\
0.499104999104999	3.1230288861673\\
0.499657999658	3.39591311370299\\
0.499938999939	3.84207182317968\\
0.499998999999	4.75342410672217\\
};
\end{axis}
\end{tikzpicture}%
}
	~	\subfloat[periodization mapping $\tau_{2,\Delta}$]{
    		\setlength\fwidth{3.5cm}
\definecolor{mycolor1}{rgb}{0.00000,0.44700,0.74100}%
\begin{tikzpicture}
\begin{axis}[%
width=0.951\fwidth,
height=\fheight,
at={(0\fwidth,0\fheight)},
scale only axis,
xmin=0,
xmax=1,
xtick={0, 0.25, 0.5, 0.75, 1},
ymin=-0.5,
ymax=0.5,
ytick={-0.5, -0.25, 0, 0.25, 0.5},
axis background/.style={fill=white},
xmajorgrids,
ymajorgrids,
legend style={legend cell align=left, align=left, draw=white!15!black}
]
\addplot [color=mycolor1]
  table[row sep=crcr]{%
0	-0.3\\
0.1000100010001	-0.4999799979998\\
0.6000600060006	0.4998799879988\\
1	-0.3\\
};
\end{axis}
\end{tikzpicture}%
}
	~
	\subfloat[combined mapping $\varphi_{j,\Delta}$]{
    		\setlength\fwidth{3.5cm}
\definecolor{mycolor1}{rgb}{0.00000,0.44700,0.74100}%
\begin{tikzpicture}
\begin{axis}[%
width=0.951\fwidth,
height=\fheight,
at={(0\fwidth,0\fheight)},
scale only axis,
xmin=0,
xmax=1,
xtick={0, 0.25, 0.5, 0.75, 1},
ymin=-4,
ymax=4,
axis background/.style={fill=white},
xmajorgrids,
ymajorgrids,
legend style={legend cell align=left, align=left, draw=white!15!black}
]
\addplot [color=mycolor1]
  table[row sep=crcr]{%
0	-0.841621233572915\\
0.00934700934700938	-0.910399307892649\\
0.0180600180600177	-0.978635872909026\\
0.0261880261880263	-1.0466784735522\\
0.0337710337710337	-1.11484729603821\\
0.0408420408420405	-1.18344751256303\\
0.0474280474280473	-1.25277442870041\\
0.0535520535520533	-1.32313109432815\\
0.0592340592340594	-1.39484109619377\\
0.0644910644910643	-1.46825215472613\\
0.0693380693380696	-1.54375489260941\\
0.0737890737890741	-1.62181154896959\\
0.0778560778560777	-1.70295928241323\\
0.0815490815490811	-1.78782854627132\\
0.084879084879085	-1.87725129262924\\
0.0878540878540877	-1.9722274396542\\
0.0904820904820909	-2.07408262358016\\
0.0927720927720932	-2.18468938651011\\
0.0947300947300951	-2.30655535074717\\
0.0963640963640966	-2.44354850246523\\
0.0976820976820978	-2.60187332159451\\
0.0986910986910985	-2.79216869600586\\
0.0994010994010992	-3.03622535884783\\
0.0998240998241	-3.38817222007596\\
0.0999880999881002	-4.06710867800757\\
0.1000001000001	-5.06895755939991\\
0.1002801002801	-3.25840568439032\\
0.100900100900101	-2.91120270121373\\
0.101846101846101	-2.67899238210228\\
0.103119103119103	-2.49837502834215\\
0.104722104722105	-2.34772583326286\\
0.106659106659107	-2.21680489855727\\
0.108935108935109	-2.09986821211158\\
0.111556111556111	-1.99333801619013\\
0.114528114528115	-1.89484863247447\\
0.117858117858118	-1.80271826979406\\
0.121554121554121	-1.71570261336787\\
0.125624125624126	-1.63286923170897\\
0.13007913007913	-1.55344641994879\\
0.134929134929135	-1.47684742856914\\
0.14018614018614	-1.40257181730664\\
0.145865145865146	-1.33017510822444\\
0.151981151981152	-1.25929275836622\\
0.158552158552158	-1.18958731447101\\
0.165600165600166	-1.12073503980804\\
0.173149173149173	-1.05244225094272\\
0.181229181229181	-0.984404380870686\\
0.18987518987519	-0.916316799033484\\
0.199132199132199	-0.847836923802342\\
0.209054209054209	-0.778597527151119\\
0.21971121971122	-0.70816164942789\\
0.231196231196231	-0.635986954732907\\
0.243637243637243	-0.561364654024658\\
0.257220257220257	-0.483302410326407\\
0.272234272234273	-0.400297934718007\\
0.289185289185289	-0.309763044083766\\
0.309133309133309	-0.206329884160702\\
0.335444335444335	-0.0730361613064447\\
0.396621396621397	0.235894658098899\\
0.416564416564417	0.340151641247481\\
0.433324433324433	0.430678344836787\\
0.448101448101448	0.51351058615652\\
0.461433461433462	0.591379527744937\\
0.473620473620474	0.665832599312747\\
0.484851484851485	0.737868996621046\\
0.495256495256495	0.80820249972124\\
0.504929504929505	0.877376854961699\\
0.513940513940514	0.945824819409002\\
0.522344522344523	1.01391793916361\\
0.53018553018553	1.08198795694036\\
0.537499537499538	1.15034488689691\\
0.544315544315545	1.21928030096836\\
0.55065955065955	1.28910439061525\\
0.556551556551557	1.36011424540001\\
0.562011562011562	1.43266445570396\\
0.567054567054567	1.50711285702036\\
0.571695571695572	1.58389457027023\\
0.575946575946576	1.66349347343554\\
0.57981857981858	1.7464909359666\\
0.583321583321584	1.83359815228638\\
0.586466586466586	1.92576555567881\\
0.589260589260589	2.02412158747113\\
0.591712591712592	2.13026631180483\\
0.593829593829594	2.24634958743683\\
0.595618595618595	2.37549104771498\\
0.597087597087597	2.52258402057826\\
0.598243598243599	2.69562804155657\\
0.599095599095599	2.90971336764274\\
0.5996545996546	3.19846743078545\\
0.5999385999386	3.6668046803759\\
0.5999995999996	4.79832231295132\\
0.600294600294601	3.24405710727704\\
0.600924600924601	2.90280200998426\\
0.601880601880602	2.67278665106336\\
0.603163603163603	2.4933497559403\\
0.604776604776605	2.34344812317476\\
0.606723606723607	2.21304668962313\\
0.609009609009609	2.09649349873021\\
0.611640611640611	1.99025868590817\\
0.614622614622615	1.89200392019746\\
0.617963617963618	1.800039053678\\
0.621670621670622	1.71316337942978\\
0.625752625752626	1.63043068872412\\
0.63021863021863	1.5511133566091\\
0.635080635080636	1.47459098394872\\
0.640350640350641	1.40036995801785\\
0.646041646041646	1.32803488056151\\
0.652170652170653	1.25719617317805\\
0.658755658755659	1.18751984130283\\
0.665817665817666	1.11869409957534\\
0.673381673381673	1.0504164331396\\
0.681477681477682	0.982383959007727\\
0.69014169014169	0.914285665802198\\
0.699417699417699	0.845788400167452\\
0.709360709360709	0.776518609213481\\
0.72004072004072	0.706040615040956\\
0.731552731552732	0.633800654191379\\
0.744024744024744	0.559091938681003\\
0.757645757645758	0.480906390648851\\
0.772708772708773	0.397722037854805\\
0.78972878972879	0.306905688746018\\
0.80979580979581	0.202938324100997\\
0.836438836438837	0.068038048960311\\
0.895773895773896	-0.231528281329838\\
0.915881915881916	-0.336528536026568\\
0.932732932732933	-0.427427104395344\\
0.947572947572947	-0.510490012777632\\
0.960952960952961	-0.588512783089551\\
0.973178973178973	-0.663072523258585\\
0.984442984442985	-0.735182981432273\\
0.994876994876995	-0.805567934914829\\
1	-0.841621233572915\\
};
\end{axis}
\end{tikzpicture}%
}
	\caption{The plots of the transformation and periodization mappings $\tau_1$ and $\tau_{2,\Delta}$ and the combined mapping $\varphi_{j,\Delta}$ with shift $\Delta = 0.1$.}
  \label{fig:Periodization_lognormal}
\end{figure}
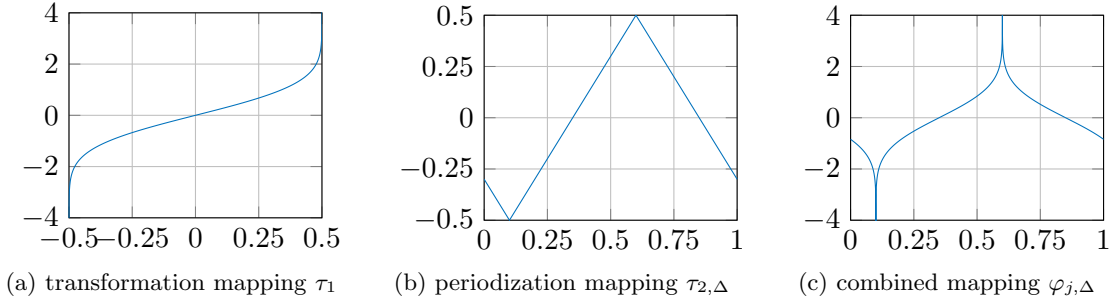

Now we define the transformation mappings $\varphi_{j,\Delta}$ for each $j=1,...,d$ for the lognormal case as
\begin{subequations}\label{eq:lognormal_transfo}
\begin{align}
&\varphi_{j,\Delta}: \T \setminus \left\lbrace \Delta,\frac12+\Delta \right\rbrace \longrightarrow \R, & \varphi_{j,\Delta}(\tilde{y}) = (\tau_1\circ\tau_{2,\Delta})(\tilde{y}),\\
&\varphi_{j,\Delta}^{-1}: \R \longrightarrow \left(\Delta,\frac12+\Delta\right), & \varphi_{j,\Delta}^{-1}(y) = (\tau_{2,\Delta}^{-1} \circ \tau_1^{-1})(y).
\end{align}
\end{subequations}
The mapping $\varphi_{j,\Delta}$ as well as its two parts $\tau_1$ and $\tau_{2,\Delta}$ are visualized in Figure \ref{fig:Periodization_lognormal}. These mappings fulfill slightly modified versions of the assumptions (A1) - (A4) taking into account the shift $\Delta$. 
Now we can use the transformation $\varphi_\Delta \coloneqq (\varphi_{j,\Delta})_{j=1}^{d}$ to receive the in $\tilde{\boldy}$ periodic signals $\tilde{u}(\boldx_g,\tilde{\boldy}) = u(\boldx_g,\varphi_\Delta(\tilde{\boldy}))$, $\boldx_g\in\mathcal{T}_G$, that can be approximated using our usFFT, cf.\ Algorithm~\ref{algo:sfft_general}. Plugging the inverse mapping $\varphi_\Delta^{-1}$ into the evaluation formula, which is similar to \eqref{eq:eval_formula}, we are now able to compute approximations of the solution functions $u_{\boldx_g}(\boldy)$ in the lognormal case as well. Again, the periodization $\varphi_{\Delta}$ is not smooth and therefore might yield non-optimal approximation results. In particular, the periodization mappings $\varphi_{j,\Delta}$ possess two poles instead of two kinks, which is a way worse smoothness behavior than in the affine case.

We now ask for the optimal choice of the parameter $\Delta$, such that the deterministic components of the sampling nodes $\tilde{y}$ of the R1Ls in the usFFT are as far as possible from $\Delta$ and $\Delta+1/2$ to reduce problems at the poles of the transformation mapping $\varphi_\Delta$. Let
\begin{align}\label{eq:lattice_nodes}
&\tilde{y}_{i,j} \coloneqq \frac{i}{M} z_j \mod 1, & i=0,...,M-1 \;\text{ and }\; j=1,...,d
\end{align}
denote the $j$-th component of the $i$-th R1L node of the $d$-dimensional R1L of size $M$. Then, we are looking for $\Delta$, such that the minimum of the two distances
\begin{align*}
\min_{\substack{i=0,...,M-1\\j=1,...,d}} \left\lvert \tilde{y}_{i,j} - \Delta  \right\rvert  \qquad \text{ and } \qquad \min_{\substack{i=0,...,M-1\\j=1,...,d}} \left\lvert \tilde{y}_{i,j} - \left( \Delta + \frac12 \right) \right\rvert
\end{align*}
is maximal.

\begin{lemma}\label{lem_delta}
Let $\Lambda(\boldz,M)$ be a $d$-dimensional R1L with prime lattice size $M\in \N$, $M>2,$ generating vector $\boldz \in \Z^d$, $z_{j_0} \not\equiv 0 \imod{M}$ for at least one $j_0\in\{1,\ldots,d\}$, and the lattice nodes $\tilde{y}_{i,j}$ as defined in \eqref{eq:lattice_nodes}. Then, we have
\begin{align*}
\Delta_{\text{opt}} \coloneqq \argmax_{0 < \Delta < \frac{1}{2M}} \left\lbrace \min \left\lbrace \min_{\substack{i=0,...,M-1\\j=1,...,d}} \left\lvert \tilde{y}_{i,j} - \Delta  \right\rvert, \min_{\substack{i=0,...,M-1\\j=1,...,d}} \left\lvert \tilde{y}_{i,j} - \left( \Delta + \frac12 \right) \right\rvert \right\rbrace \right\rbrace = \frac{1}{4M}.
\end{align*}
\end{lemma}

The proof of Lemma \ref{lem_delta} is given in Appendix \ref{app:proof_lem}.

\section{Numerics}\label{sec:numerics}

We will now test the usFFT on different, two-dimensional numerical examples. In particular, we consider the parametric PDE \eqref{eq:pde} with zero boundary condition and different random coefficients $a(\boldx,\boldy)$ and right-hand sides $f(\boldx)$.

Since our algorithm yields an approximation $u_{\boldx_g}^{\texttt{usFFT}}(\boldy)$ for each $\boldx_g \in \setT_G$ separately, we also compute the approximation error
\begin{align}
\text{err}_p^{\eta}(\boldx_g) &\coloneqq  \left( \frac{1}{n_{\text{test}}}  \sum_{j=1}^{n_{\text{test}}} \left\lvert \check{u}\left(\boldx_g,\boldy^{(j)}\right) - u^{\texttt{usFFT}}\left(\boldx_g,\boldy^{(j)}\right) \right\rvert^p \right)^{\frac{1}{p}} \\
\intertext{and}
\text{err}_\infty^{\eta}(\boldx_g) &\coloneqq  \max_{j=1,...,n_{\text{test}}} \left\lvert \check{u}\left(\boldx_g,\boldy^{(j)}\right) - u^{\texttt{usFFT}}\left(\boldx_g,\boldy^{(j)}\right) \right\rvert
\end{align}
for each $\boldx_g \in \setT_G$ separately, using $n_{\text{test}}=10^5$ different, randomly drawn test variables $\boldy^{(j)}$ from the underlying probability distribution. Here, $\check{u}(\cdot,\boldy^{(j)})$ are the finite element solutions of the PDE for fixed parameters $\boldy^{(j)}$ and $u^{\texttt{usFFT}}(\boldx_g,\cdot)$ are our approximations from the usFFT. The parameter $\eta$ denotes the used sFFT parameters as given in \autoref{table:parameters}.

\begin{table}[tb]
\caption{Parameter settings $\eta$ for the numerical tests of Algorithm \ref{algo:sfft_general}.}
\centering
\begin{tabular}{c||c||c|c||c|c||c|c|c||c|c|c||c|c} 
 $\eta$ & I & II & III & IV & V & VI & VII & VIII & IX &  X & XI & XII & XIII\\
 \hline
 $N$ & $32$ & $32$ & $64$ & $32$ & $64$ & $32$ & $64$ & $128$ & $32$ & $64$ & $128$ & $128$ & $256$ \\
 \hline
 $s, s_{\text{local}}$ & $100$ & \multicolumn{2}{c||}{$250$} & \multicolumn{2}{c||}{$500$} & \multicolumn{3}{c||}{$1000$} & \multicolumn{3}{c||}{$2000$} & \multicolumn{2}{c}{$4000$} \\
 \hline
 $\theta$ & \multicolumn{13}{c}{$1 \cdot 10^{-12}$} \\
 \hline
 $r$ & \multicolumn{13}{c}{$5$} \\
\end{tabular}
\label{table:parameters}
\end{table}

Here, $N$ is the extension of the full grid $[-N,N]^{d}$, that is used as the search space ${\Gamma \subset \Z^{d}}$. Note that we also choose $s_{\text{local}} = s$. If we miss an important frequency component at one point $\boldx_g$, it is very likely, that it is contained in the detected index set of a neighboring mesh point. Therefore, the union over all points $\boldx_g$ should be enough to avoid losing frequencies and we do not need a larger $s_{\text{local}}$. The choices for the threshold $\theta$ and the number of detection iterations $r$ are common values and the same as in \cite{KaKrVo20}. In particular, we choose the number of detection iterations $r$ as well as the probabilities $\gamma_{\operatorname{A}}$ and $\gamma_{\operatorname{B}}$ as in the case $G=1$, since we expect a huge overlap of the detected index sets and hence a small failure probability even for these parameter choices instead of the theoretical choices given in the proof of Theorem \ref{thm:complexities}. 

Further, we always use the random R1L approach in the role of Algorithm~$\operatorname{A}$ to recover the projected Fourier coefficients in the dimension-incremental method, cf.\ Section \ref{subsec:method} and \cite{KaKrVo20}. We also tested the algorithm using the single and multiple R1L approaches mentioned in Section \ref{subsec:method}, but these did not achieve significantly smaller approximation errors and are using larger numbers of samples and therefore result in longer runtimes of the algorithm. Hence, it seems reasonable to stick with the random R1L approach here. We choose the target maximum edge length of the finite element mesh $h_{\max} = 0.075$ in the FE solver. All examples consider the spatial domain $D_{\boldx}=[0,1]^2$, resulting in a finite element mesh $\setT_G \subset D_{\boldx}$ with $G=737$ inner and $104$ boundary nodes.

Further, we will also analyze the importance of and the interactions between our detected Fourier coefficients. To this end, we use the classical ANOVA decomposition of 1-periodic functions as given in \cite{PoSchmi19} or \cite{ZhaHuHoLiYa14, KuNuPlSlWa17}. Note that for instance in \cite{PoSchmi19} the ANOVA decomposition is used already in the proposed methods to receive an adaptive selection of the most important approximation terms, which we realized in our method by simply comparing the size of the projected Fourier coefficients, cf.\ Sections~\ref{subsec:method} and~\ref{subsec:usFFT}.

In particular, we consider the variance of our approximation
\begin{align*}
\sigma^2(\tilde{u}_{\boldx_g}^{\texttt{usFFT}}) \coloneqq \| \tilde{u}_{\boldx_g}^{\texttt{usFFT}} \|_{L_2(\T^{d})}^2 - \lvert c_{\boldzero}^{\texttt{usFFT}}(\tilde{u}_{\boldx_g}) \rvert^2 = \sum_{\boldk \in \I \setminus \lbrace \boldzero \rbrace} \lvert c_{\boldk}^{\texttt{usFFT}}(\tilde{u}_{\boldx_g}) \rvert^2.
\end{align*}
Now we can study different subsets $\J \subset  \I$ and estimate the variance of the approximation using only these subsets. The fraction of variance, that is explained using this subset $\J$, is then called global sensitivity index (GSI), see  \cite{So90, So01}, 
\begin{align}\label{eq:gsi}
\varrho(\J,\tilde{u}_{\boldx_g}^{\texttt{usFFT}}) \coloneqq \frac{\sigma^2(\tilde{u}_{\boldx_g, \J}^{\texttt{usFFT}})}{\sigma^2(\tilde{u}_{\boldx_g}^{\texttt{usFFT}})} = \frac{\sum_{\boldk \in \J \setminus \lbrace \boldzero \rbrace} \lvert c_{\boldk}^{\texttt{usFFT}}(\tilde{u}_{\boldx_g}) \rvert^2}{\sum_{\boldk \in \I \setminus \lbrace \boldzero \rbrace} \lvert c_{\boldk}^{\texttt{usFFT}}(\tilde{u}_{\boldx_g}) \rvert^2} \; \in \, [0,1],
\end{align}
where we define $\tilde{u}_{\boldx_g, \J}^{\texttt{usFFT}}(\boldy)\coloneqq\sum_{\boldk\in \J}c_{\boldk}^{\texttt{usFFT}}(\tilde{u}_{\boldx_g})\,\e^{2\pi\ii\boldk\cdot \varphi^{-1}(\boldy)}$.
In our examples, we will mainly consider the subsets $\J_{\ell}$ of all frequencies $\boldk$ with exactly $\ell$ non-zero components, i.e.,
\begin{align}\label{eq:J}
\J_{\ell} \coloneqq \left\lbrace \boldk \in \I: \|\boldk\|_0 \coloneqq \lvert \lbrace i \in \lbrace 1,...,d \rbrace: k_i \not = 0\rbrace\rvert = \ell \right\rbrace,
\end{align}
but of course several other choices of $\J_{\ell}$ might be interesting as well for different applications.

Finally, one can also think about evaluating various quantities of interest of the approximation. Here, we will consider the expectation value $\E(u_{\boldx_g}^{\texttt{usFFT}})$ as one example of such quantities. We use a Monte-Carlo approximation of the expectation value
\begin{align*}
\overline{\check{u}_{\boldx_g}} \coloneqq \frac{1}{n_{\text{MC}}} \sum_{j=1}^{n_{\text{MC}}} \check{u}\left(\boldx_g,\boldy^{(j)}\right)
\end{align*}
of the finite element approximation using ${n_{\text{MC}}}$ random samples for comparison.

\subsection{Expectation value of the approximation}
\label{subsec:moments}

Computing the expectation value of our approximation $u_{\boldx_g}^{\texttt{usFFT}}$ requires some additional effort, depending on the particular model and eventually used periodization methods. By definition, the expectation value is given by
\begin{align*}
\E(u_{\boldx_g}^{\texttt{usFFT}}) \coloneqq \int_{D_{\boldy}} u_{\boldx_g}^{\texttt{usFFT}}(\boldy) \d\mu(\boldy) = \int_{D_{\boldy}} u_{\boldx_g}^{\texttt{usFFT}}(\boldy) \,p(\boldy) \d\boldy,
\end{align*}
where $p$ is the probability density function of the random variable $\boldy$.

For the periodic model, we do not need any periodization. Therefore, the approximation of the solution reads as
\begin{align*}
u_{\boldx_g}^{\texttt{usFFT}}(\boldy) = \sum_{\boldk \in \I} c_{\boldk}^{\texttt{usFFT}}(u_{\boldx_g})\,\e^{2\pi\ii\boldk\cdot\boldy}
\end{align*}
with $\I$ the frequency set and $c_{\boldk}^{\texttt{usFFT}}(u_{\boldx_g})$ the corresponding approximated Fourier coefficients computed by the usFFT.
The random variable $\boldy$ is assumed to be uniformly distributed in $D_{\boldy} = [-\frac12,\frac12]^{d}$ in this case. Hence, we have
\begin{align*}
\E(u_{\boldx_g}^{\texttt{usFFT}}) &= \int_{D_{\boldy}} u_{\boldx_g}^{\texttt{usFFT}}(\boldy) \,p(\boldy) \d\boldy \\
&= \int_{[-\frac12,\frac12]^{d}} 1^{-d} \, \sum_{\boldk \in \I} c_{\boldk}^{\texttt{usFFT}}(u_{\boldx_g})\,\e^{2\pi\ii\boldk\cdot\boldy} \d\boldy \\
&= \sum_{\boldk \in \I} c_{\boldk}^{\texttt{usFFT}}(u_{\boldx_g}) \int_{[-\frac12,\frac12]^{d}} \e^{2\pi\ii\boldk\cdot\boldy} \d\boldy \\
&= \sum_{\boldk \in \I} c_{\boldk}^{\texttt{usFFT}}(u_{\boldx_g}) \, \delta_{\boldk} = c_{\boldzero}^{\texttt{usFFT}}(u_{\boldx_g}).
\end{align*}

In the affine case, we use the tent transformation \eqref{eq:tent}, such that our approximation reads as
\begin{align*}
u_{\boldx_g}^{\texttt{usFFT}}(\boldy) = \sum_{\boldk \in \I} c_{\boldk}^{\texttt{usFFT}}(\tilde{u}_{\boldx_g})\,\e^{\pi\ii\boldk\cdot\frac{\boldy - \alpha \boldone}{\beta-\alpha}}
\end{align*}
with $\boldone = (1,1,...,1) \in \R^{d}$. Again, the random variable $\boldy$ is assumed to be uniformly distributed, but for this computation we work with the more general domain $D_{\boldy} = [\alpha,\beta]^{d}$. Therefore, we have
\begin{align}
\E(u_{\boldx_g}^{\texttt{usFFT}}) &= \int_{D_{\boldy}} u_{\boldx_g}^{\texttt{usFFT}}(\boldy) \,p(\boldy) \d\boldy \nonumber \\
&= \int_{[\alpha,\beta]^{d}} (\beta-\alpha)^{-d} \sum_{\boldk \in \I} c_{\boldk}^{\texttt{usFFT}}(\tilde{u}_{\boldx_g})\,\e^{\pi\ii\boldk\cdot\frac{\boldy - \alpha \boldone}{\beta-\alpha}} \d\boldy \nonumber \\
&= \sum_{\boldk \in \I} c_{\boldk}^{\texttt{usFFT}}(\tilde{u}_{\boldx_g}) (\beta-\alpha)^{-d} \int_{[\alpha,\beta]^{d}} \e^{\pi\ii\boldk\cdot\frac{\boldy - \alpha \boldone}{\beta-\alpha}} \d\boldy \nonumber \\
&= \sum_{\boldk \in \I} c_{\boldk}^{\texttt{usFFT}}(\tilde{u}_{\boldx_g}) \, D_{\boldk} \label{eq:affine_ew}
\end{align}
with
\begin{align*}
D_{\boldk} \coloneqq \prod_{j=1}^{d} D_{k_j} \qquad \text{and} \qquad D_{k_j} \coloneqq \begin{cases}
\frac{2\ii}{\pi k_j} & k_j \equiv 1 \mod 2\\
1 & k_j = 0 \\
0 & \text{else}.
\end{cases}
\end{align*}
Note that the parameters $\alpha$ and $\beta$ vanish completely. Thus, the formula is independent of the particular domain $D_{\boldy} = [\alpha,\beta]^{d}$.

Finally, the lognormal model involves the more complicated transformation mappings $\varphi_{j,\Delta}$ given in \eqref{eq:lognormal_transfo}. Thus, the approximation reads as
\begin{align*}
u_{\boldx_g}^{\texttt{usFFT}}(\boldy) = \sum_{\boldk \in \I} c_{\boldk}^{\texttt{usFFT}}(\tilde{u}_{\boldx_g})\,\e^{2\pi\ii\boldk\cdot\varphi_{\Delta}^{-1}(\boldy)}.
\end{align*}
Here, the random variable $\boldy$ is standard normally distributed, i.e., $\boldy \sim \mathcal{N}(\boldzero,\boldI)$ with $\boldI$ the identity matrix of dimension $d$. Hence, the expectation value can be written as
\begin{align*}
\E(u_{\boldx_g}^{\texttt{usFFT}}) &= \int_{D_{\boldy}} u_{\boldx_g}^{\texttt{usFFT}}(\boldy) \,p(\boldy) \d\boldy \\
&= \int_{\R^{d}} (2\pi)^{-\frac{d}{2}} \, \e^{-\frac12 \|\boldy\|^2} \sum_{\boldk \in \I} c_{\boldk}^{\texttt{usFFT}}(\tilde{u}_{\boldx_g})\,\e^{2\pi\ii\boldk\cdot\varphi_{\Delta}^{-1}(\boldy)} \d\boldy \\
&= \sum_{\boldk \in \I} c_{\boldk}^{\texttt{usFFT}}(\tilde{u}_{\boldx_g}) (2\pi)^{-\frac{d}{2}} \int_{\R^{d}} \e^{-\frac12 \|\boldy\|^2} \, \e^{2\pi\ii\boldk\cdot\varphi_{\Delta}^{-1}(\boldy)} \d\boldy \\
&= \sum_{\boldk \in \I} c_{\boldk}^{\texttt{usFFT}}(\tilde{u}_{\boldx_g}) D_{\boldk,\Delta} \\
\end{align*}
with
\begin{align*}
D_{\boldk,\Delta} \coloneqq \prod_{j=1}^{d} D_{k_j,\Delta} \qquad \text{and} \qquad D_{k_j,\Delta} \coloneqq \begin{cases}
\frac{2 \ii}{\pi k_j} \e^{2\pi\ii k_j \Delta}& k_j \equiv 1 \mod 2\\
1 & k_j = 0 \\
0 & \text{else}.
\end{cases}
\end{align*}
Note that the factors $D_{k_j,\Delta}$ are exactly the same as the $D_{k_j}$ in the affine case up to the correction term $\e^{2\pi\ii k_j \Delta}$ due to the shift with $\Delta$.

\subsection{Periodic example} \label{subsec:periodic}

We consider the example from \cite[Sec.~6]{KaKaKu20} using the domain $D_{\boldx} = (0,1)^2$ with right-hand side $f(\boldx) = x_2$ and the random coefficient
\begin{align*}
& a(\boldx,\boldy) \coloneqq 1 + \frac{1}{\sqrt{6}} \sum_{j=1}^{d} \sin(2\pi y_j) \,\psi_j(\boldx), & \boldx \in D_{\boldx}, \, \boldy \in D_{\boldy},
\end{align*}
with the random variables $\boldy \sim \mathcal{U}\left([-\frac12,\frac12]^{d}\right)$ and
\begin{align*}
& \psi_j(\boldx) \coloneqq c j^{-\mu} \sin(j\pi x_1) \sin(j\pi x_2), & \boldx \in D_{\boldx},\, j\geq 1,
\end{align*}
where $c>0$ is a constant and $\mu>1$ is the decay rate.  Accordingly, we get
\begin{align*}
a_{\min} = 1-\frac{c}{\sqrt6}\zeta(\mu) \qquad \text{and} \qquad a_{\max} = 1+\frac{c}{\sqrt6}\zeta(\mu),
\end{align*}
such that for $c < \frac{\sqrt6}{\zeta(\mu)}$ the uniform ellipticity assumption \eqref{eq:UEA} is fulfilled.

We test the usFFT with the stochastic dimension $d = 10$ on the two parameter choices $\mu = 1.2,\; c = 0.4$ and $\mu = 3.6,\; c = 1.5$ from \cite{KaKaKu20}. The first choice seems to model a more difficult PDE, since the decay of the functions $\psi_j$ w.r.t.\ $j$ is very slow and we have $a_{\min} = 0.08690$ and $a_{\max} = 1.91310$. This range of $a$ is wider and $a_{\min}$ is closer to zero than for the quickly decaying second parameter choice with $a_{\min} = 0.31660$ and $a_{\max} = 1.68340$. 

\begin{figure}[tb]
	\centering	
	\subfloat[$\overline{\check{u}_{\boldx_g}}$]{
		\includegraphics[scale=0.35]{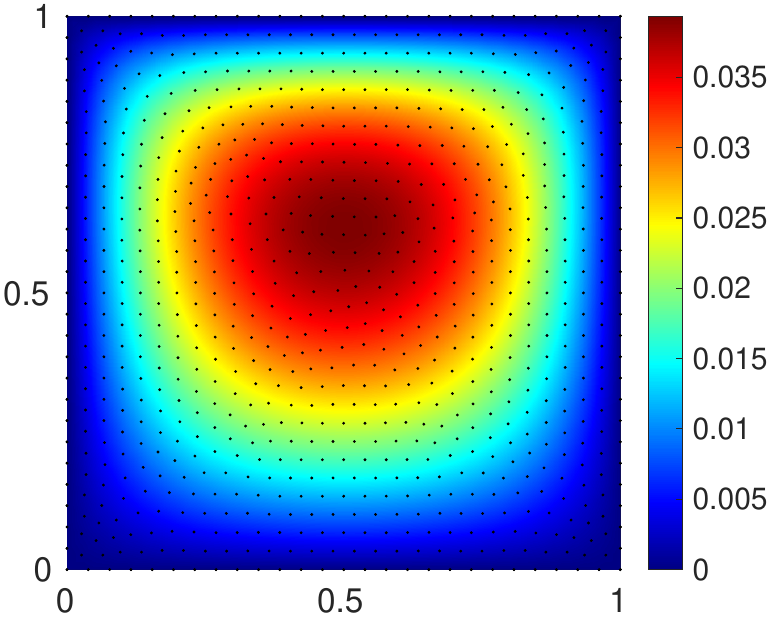} 
	}
	~	
	\subfloat[$\err_1^{\eta}(\boldx_g)$]{
    		\includegraphics[scale=0.35]{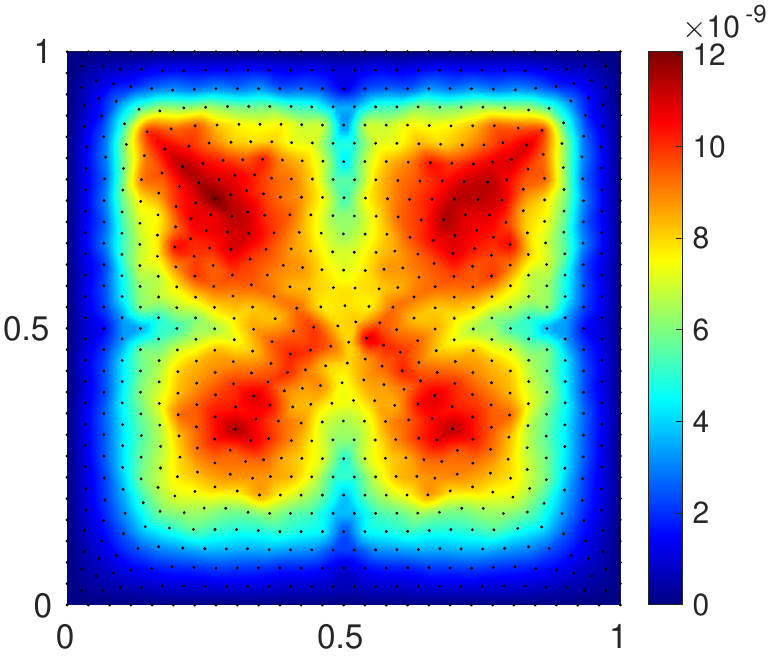} 
	}
	~
	\subfloat[$\err_2^{\eta}(\boldx_g)$]{
    		\includegraphics[scale=0.35]{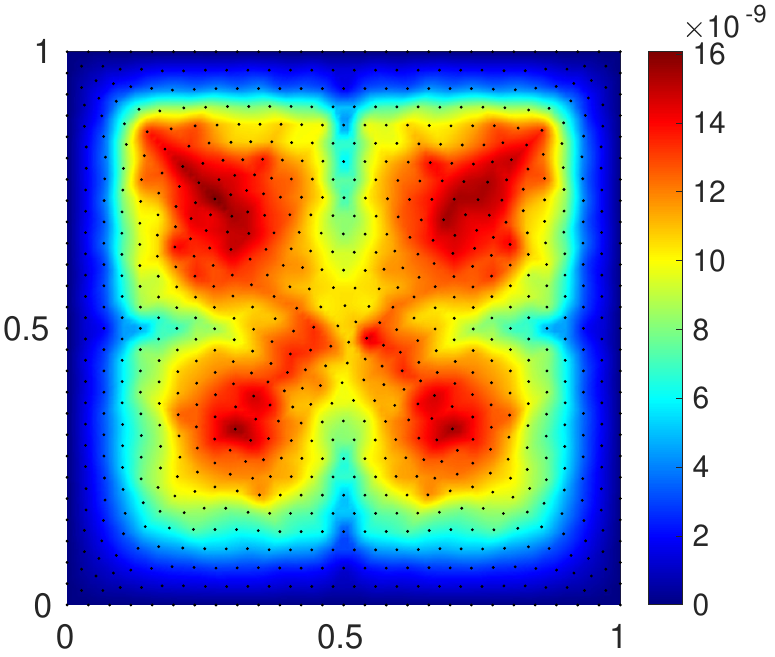} 
	}
	\caption{The MC approximation $\overline{\check{u}_{\boldx_g}}$ and the approximation errors $\err_1^{\eta}(\boldx_g)$ and $\err_2^{\eta}(\boldx_g)$ for the periodic example with $\mu = 1.2$, $c=0.4$, $d = 10$, $\eta=\textnormal{VII}$, i.e., $s=1000$, $N=64$.}
  \label{fig:periodic_example}
\end{figure}

\begin{figure}[tb]
	\centering	
    	\includegraphics[scale=0.65]{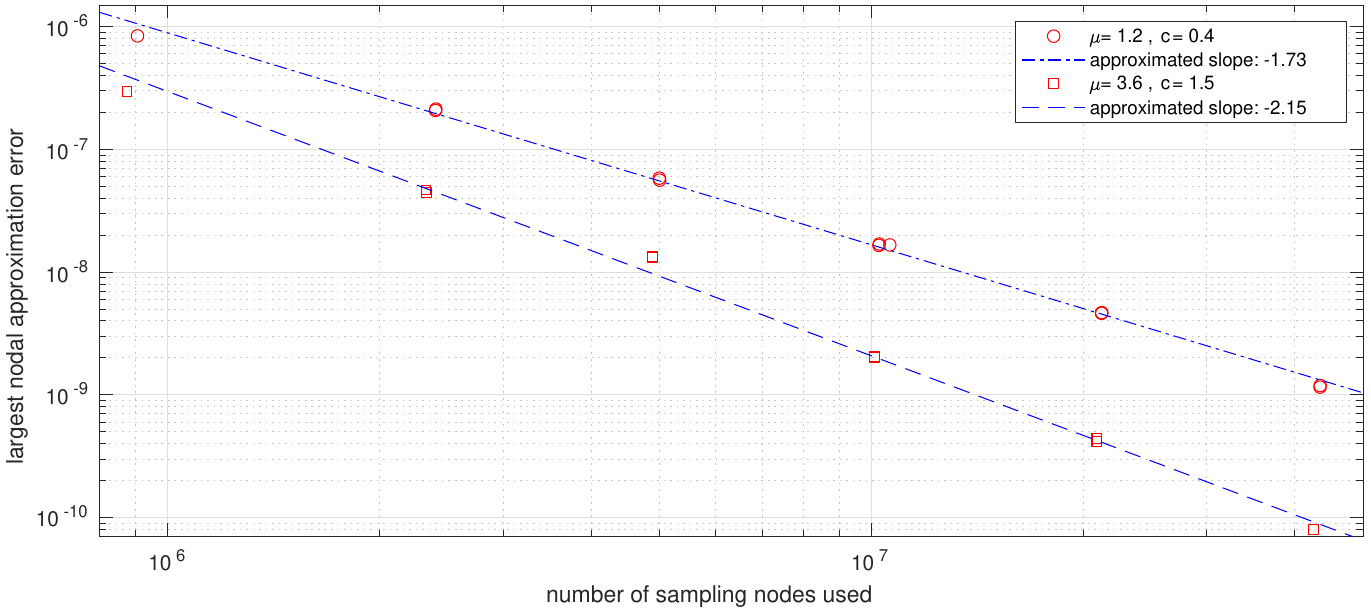} 
	\caption{Largest error $\err_2^{\eta}$ w.r.t.\ the nodes $\boldx_g$ for all parameter settings $\eta$ displayed in Table \ref{table:parameters} for the periodic example.}
  \label{fig:periodic_decay}
\end{figure}

Figure \ref{fig:periodic_example} illustrates the total approximation error $\err_p^{\eta}(\boldx_g)$ for $p=1$ and $p=2$ as well as the Monte-Carlo approximation of the expectation value $\overline{\check{u}_{\boldx_g}}$ using $n_{\text{MC}}=10^6$ samples for comparison.
A more detailed insight on the decay of the error is given in Figure \ref{fig:periodic_decay}. There, the largest approximation error $\err_2^{\eta}$ w.r.t.\ the nodes $\boldx_g$ is given with the number of samples used in the corresponding usFFT. Note that this number scales directly with the sparsity parameter $s$, while the extension parameter $N$ has nearly no impact. Hence, the data points in Figure \ref{fig:periodic_decay} are ordered from left to right from $s=100$ to $s=4000$.
Finally, Figure \ref{fig:periodic_anova} shows the cardinality of the sets $\J_\ell$, i.e., the number of frequencies detected with exactly $\ell$ non-zero components as given in \eqref{eq:J}, as well as the corresponding global sensitivity indices $\varrho(\J_\ell,u_{\boldx_g}^{\texttt{usFFT}})$ given in \eqref{eq:gsi}. Note that these values also depend on the considered point $\boldx_g$. Therefore, the bars show the smallest and largest GSI among all nodes $\boldx_g \in \setT_G$ as well as their median and mean value.

\begin{figure}[tb]
	\centering	
	\setlength\fwidth{12cm}
	\setlength\fheight{7cm}
\definecolor{mycolor1}{rgb}{0.00000,0.44700,0.74100}%
\definecolor{mycolor2}{rgb}{0.85000,0.32500,0.09800}%
\definecolor{mycolor3}{rgb}{0.92900,0.69400,0.12500}%
\begin{tikzpicture}
\pgfplotsset{
	scale only axis,
	xmin = 0.5,
	xmax = 4.5,
	xtick={1, 2, 3, 4},
	every x tick label/.append style={font=\footnotesize},
	xlabel style={font=\footnotesize\color{white!15!black}},
	xlabel={number $\ell$ of non-zero frequency components},
	xmajorgrids,
	width=\fwidth,
	height=\fheight,
	at={(0\fwidth,0\fheight)},
	legend style={legend cell align=left, align=left, draw=white!15!black}
}

\begin{axis}[
set layers,
axis y line*=left,
every outer y axis line/.append style={mycolor1},
every y tick label/.append style={font=\footnotesize\color{mycolor1}},
every y tick/.append style={mycolor1},
ymode=linear,
ymin=0,
ymax=2500,
ylabel style={font=\footnotesize\color{mycolor1}},
ylabel={number $\vert\J_{\ell}\vert$ of detected frequencies},
yticklabel pos=left,
]

\addplot[ybar, bar width=0.35, fill=mycolor1, draw=black, area legend] table[row sep=crcr] {%
0.8	68\\
1.8	820\\
2.8	2317\\
3.8	1614\\
};
\node[above, align=center]
at (axis cs:0.8,68) {\small 68};
\node[above, align=center]
at (axis cs:1.8,820) {\small 820};
\node[above, align=center]
at (axis cs:2.8,2317) {\small 2317};
\node[above, align=center]
at (axis cs:3.8,1614) {\small 1614};

\end{axis}

\begin{axis}[%
set layers,
log origin=infty,
clip = false,
axis y line*=right,
every outer y axis line/.append style={mycolor2},
every y tick label/.append style={font=\footnotesize\color{mycolor2}},
every y tick/.append style={mycolor2},
ymode=log,
ymin=1e-10,
ymax=1,
ytick={1e-10, 1e-9, 1e-8, 1e-7, 1e-6, 1e-5, 1e-4, 1e-3, 1e-2, 1e-1, 1},
yminorticks=true,
ylabel style={font=\footnotesize\color{mycolor2}},
ylabel={global sensitivity indices $\varrho(\J_{\ell},\hat{u}_{\boldx_g})$},
yticklabel pos=right,
ymajorgrids,
]

\addplot[ybar, bar width=0.35, fill=mycolor2, draw=black, area legend, postaction={line space=12pt, pattern=my north east lines}, forget plot] table[row sep=crcr] {%
1.2	0.999988154256956\\
2.2	0.00764036471255735\\
3.2	1.27220384916877e-05\\
4.2	3.13192860574541e-08\\
};
\node[above, align=center]
at (axis cs:1.2,1) {\tiny 0.9999};
\node[above, align=center]
at (axis cs:2.2,0.008) {\tiny 0.0076};
\node[above, align=center]
at (axis cs:3.2,1.2e-05) {\tiny 1.27e-05};
\node[above, align=center]
at (axis cs:4.2,3.1e-8) {\tiny 3.13e-08};

\addplot[ybar, bar width=0.35, fill=mycolor3, draw=black, area legend, forget plot] table[row sep=crcr] {%
1.2	0.992353812566311\\
2.2	1.17471273553332e-05\\
3.2	4.09113257633601e-08\\
4.2	1.05161349819723e-10\\
};
\node[below, align=center]
at (axis cs:1.2,0.992) {\tiny 0.9923};
\node[below, align=center]
at (axis cs:2.2,1.1e-05) {\tiny 1.17e-05};
\node[below, align=center]
at (axis cs:3.2,4.1e-08) {\tiny 4.09e-08};
\node[below, align=center]
at (axis cs:4.2,1e-10) {\tiny 1.05e-10};

\addplot [color=black, dashed, forget plot]
  table[row sep=crcr]{%
1.025	0.998832323658039\\
1.375	0.998832323658039\\
};
\addplot [color=black, dotted, forget plot]
  table[row sep=crcr]{%
1.025	0.999272250936575\\
1.375	0.999272250936575\\
};

\addplot [color=black, dashed, line width=0.5mm]
  table[row sep=crcr]{%
2.025	0.00116576877836653\\
2.375	0.00116576877836653\\
};
\addplot [color=black, dotted, line width=0.5mm]
  table[row sep=crcr]{%
2.025	0.000727087112337809\\
2.375	0.000727087112337809\\
};
\legend{mean, median}
\addplot [color=black, dashed, line width=0.5mm, forget plot]
  table[row sep=crcr]{%
3.025	1.90448499093387e-06\\
3.375	1.90448499093387e-06\\
};
\addplot [color=black, dotted, line width=0.5mm, forget plot]
  table[row sep=crcr]{%
3.025	9.17434242444728e-07\\
3.375	9.17434242444728e-07\\
};
\addplot [color=black, dashed, line width=0.5mm, forget plot]
  table[row sep=crcr]{%
4.025	3.07860410449744e-09\\
4.375	3.07860410449744e-09\\
};
\addplot [color=black, dotted, line width=0.5mm, forget plot]
  table[row sep=crcr]{%
4.025	1.17528203156181e-09\\
4.375	1.17528203156181e-09\\
};

\end{axis}
\end{tikzpicture}%
\caption{Cardinality (left, blue, solid) of the index sets $\J_\ell$ and the corresponding largest (right, orange, striped), smallest (right, yellow, solid), mean (dashed line) and median (dotted line) of the global sensitivity indices $\varrho(\J_\ell,u_{\boldx_g}^{\texttt{usFFT}})$ w.r.t.\ $\boldx_g$ for the periodic example with $\mu = 1.2$, $c=0.4$, $d = 10$, $\eta=\textnormal{XI}$, i.e., $s=2000$, $N=128$.}
  \label{fig:periodic_anova}
\end{figure}

\subsubsection*{Discussion}

The absolute error $\err_p^{\eta}$ in Figure \ref{fig:periodic_example} is very small compared to the function values of $\overline{\check{u}_{\boldx_g}}$. Thus, our approximation $u_{\boldx_g}^{\texttt{usFFT}}$ is already a very good approximation for these relatively small sparsity parameters $s$ and extension parameters $N$.  

The periodic setting results in very quickly decaying Fourier coefficients. Obviously, the same holds for their projections computed in the dimension-incremental steps. In particular, most of the one-dimensional projections in step~1 {\&}~2a of Algorithm~\ref{algo:sfft_general}, e.g., all projections with component $k_t$ with $\lvert k_t \rvert >4$, at the start of each iteration are so small, that they are neglected immediately. Hence, the one-dimensional index sets $\I^{(t)}$ are independent of $N$ (for large enough $N$) and so the choice of $N$ has only a marginal impact. Note that we also tested our algorithm with smaller thresholds $\theta$, but the additionally detected and not neglected frequencies did not change the approximation significantly in the end.

We indicate some kind of linear behavior in the double logarithmic Figure \ref{fig:periodic_decay}. Additional tests showed, that even for smaller sparsity parameters $1\leq s<100$ the corresponding samples-error-pair fits into this model, i.e., there seems to be no pre-asymptotic behavior of our algorithm. In \cite{KaKaKu20} the theoretical decay rates are often smaller than the error decay observed in numerical experiments. We also observe a relatively fast decay compared to these theoretical rates. On the other hand, the decay of the approximation error $\err_2^{\eta}$ for the faster decaying random coefficient $a$ with $\mu = 3.6$ is not that much better than the decay of the more complicated example with $\mu = 1.2$. It seems like our algorithm is capable of handling the more difficult problem very well, but also does not yield that much further advantages when being applied to easier problems, i.e., with larger $\mu$, larger $a_{\min}$ and a smaller range of the interval $[a_{\min}, a_{\max}]$. Note that most of the samples needed are required for the detection of the frequency set $\I$ and only a small fraction is really used for the final computation of the corresponding Fourier coefficients, cf.\ Section \ref{subsec:known_freq} and Remark \ref{rem:sizes}. We also computed the approximation error $\text{err}_\infty^{\eta}$ for different $\eta$ for both parameter choices of $\mu$ and $c$. Obviously, these errors have to be larger than the shown errors $\err_2^{\eta}$, but the actual magnitude of $\text{err}_\infty^{\eta}$ is only about $10$ or $15$ times as large as the errors $\err_2^{\eta}$. Hence, the pointwise approximation error seems to stay in a reasonable size for any randomly drawn $\boldy$.

As we saw in Section \ref{subsec:moments}, the expectation value of our approximation $\E(u_{\boldx_g}^{\texttt{usFFT}})$ is simply its zeroth Fourier coefficient $c_{\boldzero}^{\texttt{usFFT}}(u_{\boldx_g})$. Since this coefficient is included and computed for each sparsity parameter $s$ anyway, it seems like our different parameter choices would not influence the precision of its approximation at first sight. But for larger sparsity parameters $s$, we compute more Fourier coefficients in our algorithm, where possible aliasing effects should spread evenly among all of these coefficients, i.e., the particular so-called aliasing error on $c_{\boldzero}^{\texttt{usFFT}}(u_{\boldx_g})$, cf.\ \cite{KaPoVo17, KaKrVo20}, gets smaller and the approximation improves. Unfortunately, this is not visible in our numerical tests, since the comparison value $\overline{\check{u}_{\boldx_g}}$ behaves too poorly. In detail, we would have to investigate very small sparsity parameters $s<25$ to observe the described effects. For all of our parameter choices $\eta$, the Monte-Carlo approximation $\overline{\check{u}_{\boldx_g}}$ with ${n_{\text{MC}}} = 5\cdot 10^6$ samples is not accurate enough to give insight on the particular behavior of our approximation of the expectation value.

Figure \ref{fig:periodic_anova} shows, that there are no frequencies detected with all or nearly all components being active. Further, even though only $68$ of the $4819$ frequencies detected (excluding $c_{\boldzero}$) have exactly one non-zero component, i.e., are supported on the axis cross, they contain more than $99\%$ of the variance of our approximation. So the higher-dimensional frequencies with two, three or four non-zero components seem to be nearly neglectable for the approximation.

\subsection{Affine example} \label{subsec:affine}

For the affine case, we consider an example from \cite[Sec.~11]{EiGiSchwZa14} with domain $D_{\boldx} = (0,1)^2$, right-hand side $f(\boldx) \equiv 1$ and the random coefficient
\begin{align*}
&a(\boldx,\boldy) \coloneqq 1 + \sum_{j=1}^{d} y_j \psi_j(\boldx), & \boldx \in D_{\boldx}, \, \boldy \in D_{\boldy},
\end{align*}
with the random variables $\boldy \sim \mathcal{U}([-1,1]^{d})$ and
\begin{align*}
&\psi_j(\boldx) \coloneqq c j^{-\mu} \cos(2\pi m_1(j)\,x_1)\, \cos(2\pi m_2(j)\,x_2), & \boldx \in D_{\boldx},\, j\geq 1,
\end{align*}
where again $c > 0$ is a constant and $\mu > 1$ the decay rate. Further, $m_1(j)$ and $m_2(j)$ are defined as
\begin{align*}
m_1(j) \coloneqq j-\frac{k(j) (k(j)+1)}{2} \quad \text{and} \quad m_2(j) \coloneqq k(j)-m_1(j)\\
\end{align*}
with $k(j) \coloneqq \lfloor -1/2 + \sqrt{1/4 + 2j} \rfloor$. Table \ref{table:diagonal} shows the numbers $m_1(j), m_2(j)$ and $k(j)$ for a few $j \geq 1$. 

\begin{table}[tb]
\caption{The values of $m_1(j), m_2(j)$ and $k(j)$.}
\centering
\begin{tabular}{c||c|c||c|c|c||c|c|c|c||c|c|c|c|c||c} 
 $j$ & $1$ & $2$ & $3$ & $4$ & $5$ & $6$ & $7$ & $8$ & $9$ & $10$ & $11$ & $12$ & $13$ & $14$ & \ldots \\
 \hline\hline
 $m_1(j)$ & $0$ & $1$ & $0$ & $1$ & $2$ & $0$ & $1$ & $2$ & $3$ & $0$ & $1$ & $2$ & $3$ & $4$ & \ldots \\
 \hline
 $m_2(j)$ & $1$ & $0$ & $2$ & $1$ & $0$ & $3$ & $2$ & $1$ & $0$ & $4$ & $3$ & $2$ & $1$ & $0$ & \ldots \\
 \hline
 $k(j)$ & \multicolumn{2}{c||}{$1$} & \multicolumn{3}{c||}{$2$} & \multicolumn{4}{c||}{$3$} & \multicolumn{5}{c||}{$4$} & \ldots \\
\end{tabular}
\label{table:diagonal}
\end{table}

As before, we get that $a_{\min} = 1-c\,\zeta(\mu)$ and $a_{\max} = 1+c\,\zeta(\mu)$, such that for $c < \frac{1}{\zeta(\mu)}$ the uniform ellipticity assumption \eqref{eq:UEA} is fulfilled. Here, we use the parameter choices from \cite{EiGiSchwZa14} with $\mu = 2$ for a relatively slow decay and $c = \frac{0.9}{\zeta(2)} \approx 0.547$ to end up with $a_{\min}=0.1$ and $a_{\max}=1.9$, which is very similar to the first parameter choice in the periodic case. We choose the stochastic dimension $d=20$ as in \cite{EiGiSchwZa14}.

\begin{figure}[tb]
	\centering	
	\subfloat[$\overline{\check{u}_{\boldx_g}}$]{
		\includegraphics[scale=0.35]{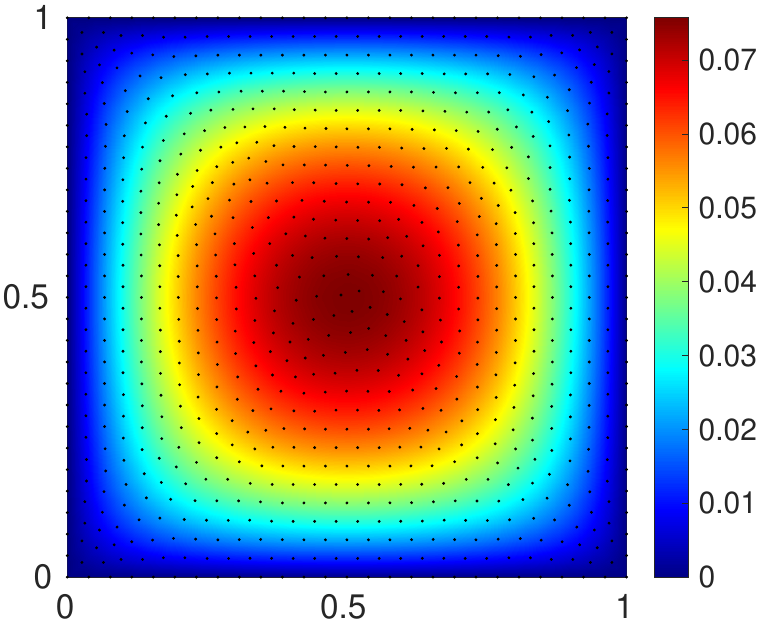} \label{fig:affine_ew1}
	}
	~	
	\subfloat[$\eta =$ I ($s=100, N=32$)]{
    		\includegraphics[scale=0.35]{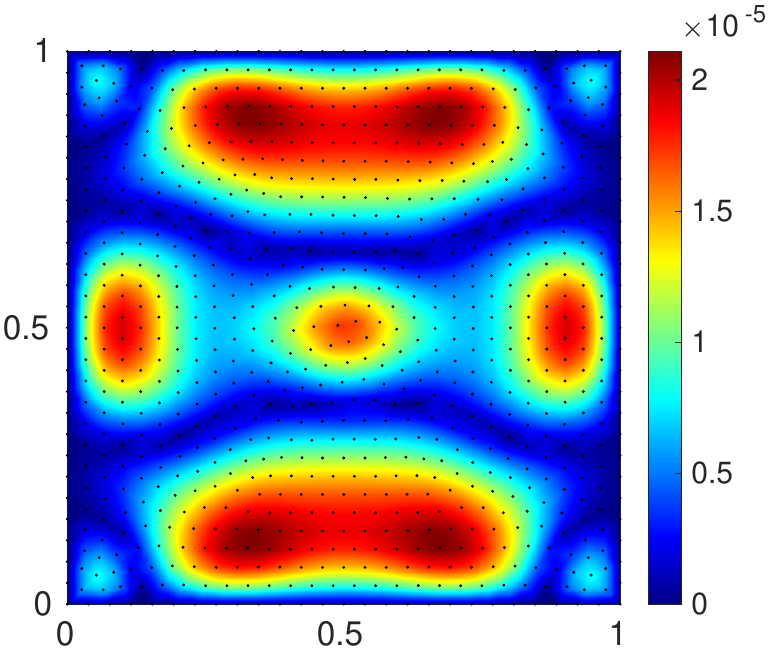} 
	}
	~
	\subfloat[$\eta =$ II ($s=250, N=32$)]{
    		\includegraphics[scale=0.35]{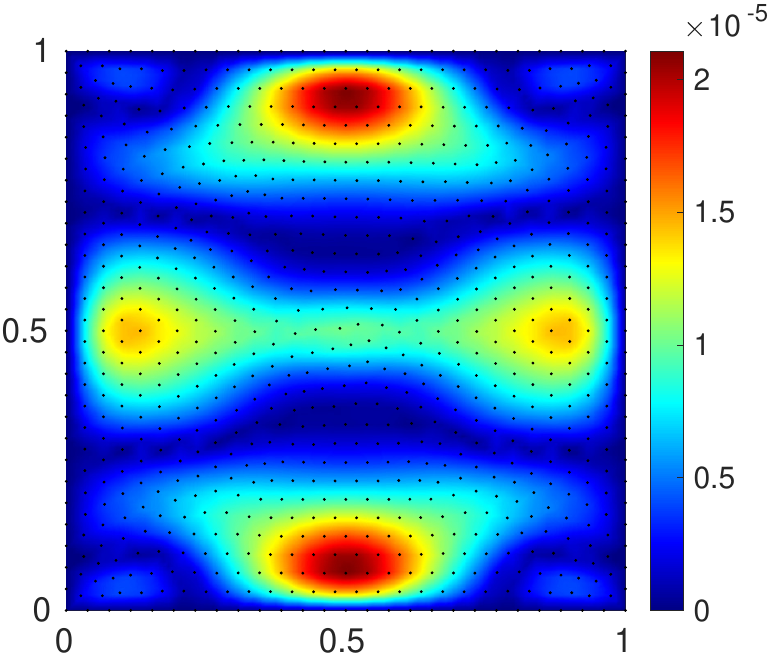} 
	}
	\caption{The MC approximation $\overline{\check{u}_{\boldx_g}}$ and the pointwise errors $\vert \overline{\check{u}_{\boldx_g}} - \E(u_{\boldx_g}^{\texttt{usFFT}})\vert$ for $\eta =$ I and II for the affine example.}
  \label{fig:affine_ew}
\end{figure}

Figure \ref{fig:affine_ew} illustrates the Monte-Carlo approximation of the expectation value $\overline{\check{u}_{\boldx_g}}$ with ${n_{\text{MC}}} = 10^6$ samples used as well as the pointwise error $\vert \overline{\check{u}_{\boldx_g}} - \E(u_{\boldx_g}^{\texttt{usFFT}})\vert$ for two different parameter choices $\eta$ with $\E(u_{\boldx_g}^{\texttt{usFFT}})$ as given in \eqref{eq:affine_ew}. 
\begin{figure}[tb]
	\centering	
    	\includegraphics[scale=0.65]{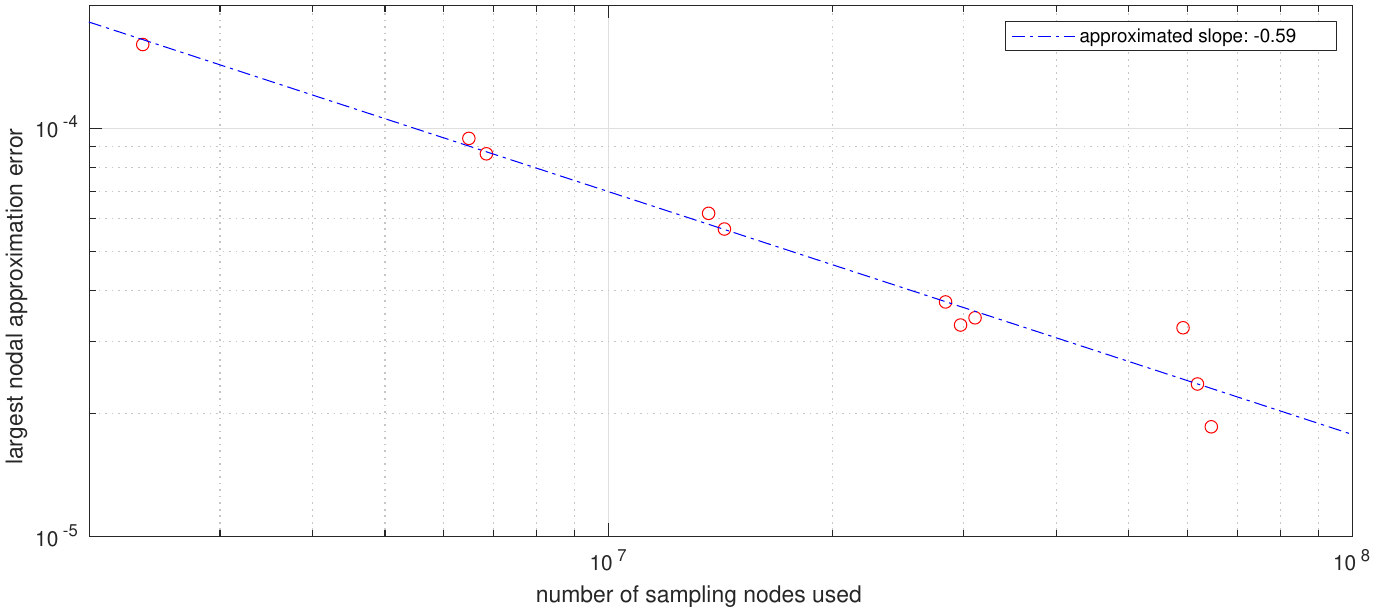} 
	\caption{Largest error $\err_2^{\eta}$ w.r.t.\ the nodes $\boldx_g$ for the parameter settings $\eta=$ I to XI displayed in Table \ref{table:parameters} for the affine example.}
  \label{fig:affine_decay}
\end{figure}
Figure \ref{fig:affine_decay} again shows the largest error $\err_2^{\eta}$ w.r.t.\ the nodes $\boldx_g$ for different parameter settings $\eta$. This time, we can observe a small increase in the number of used samples for larger extensions $N$, which was not visible in the periodic example. Hence, the parameter settings $\eta = $ I to XI have monotonously increasing sampling sizes, i.e., the data points in Figure \ref{fig:affine_decay} are ordered from left to right w.r.t.\ increasing $\eta$ this time.
\begin{figure}[tb]
	\centering	
	\setlength\fwidth{12cm}
	\setlength\fheight{7cm}
\definecolor{mycolor1}{rgb}{0.00000,0.44700,0.74100}%
\definecolor{mycolor2}{rgb}{0.85000,0.32500,0.09800}%
\definecolor{mycolor3}{rgb}{0.92900,0.69400,0.12500}%
\begin{tikzpicture}
\pgfplotsset{
	scale only axis,
	xmin = 0.5,
	xmax = 4.5,
	xtick={1, 2, 3, 4},
	every x tick label/.append style={font=\footnotesize},
	xlabel style={font=\footnotesize\color{white!15!black}},
	xlabel={number $\ell$ of non-zero frequency components},
	xmajorgrids,
	width=\fwidth,
	height=\fheight,
	at={(0\fwidth,0\fheight)},
	legend style={legend cell align=left, align=left, draw=white!15!black}
}

\begin{axis}[
set layers,
axis y line*=left,
every outer y axis line/.append style={mycolor1},
every y tick label/.append style={font=\footnotesize\color{mycolor1}},
every y tick/.append style={mycolor1},
ymode=linear,
ymin=0,
ymax=3000,
ylabel style={font=\footnotesize\color{mycolor1}},
ylabel={number $\vert\J_{\ell}\vert$ of detected frequencies},
yticklabel pos=left,
]

\addplot[ybar, bar width=0.35, fill=mycolor1, draw=black, area legend] table[row sep=crcr] {%
0.8	470\\
1.8	2568\\
2.8	1246\\
3.8	54\\
};
\node[above, align=center]
at (axis cs:0.8,470) {\small 470};
\node[above, align=center]
at (axis cs:1.8,2568) {\small 2568};
\node[above, align=center]
at (axis cs:2.8,1246) {\small 1246};
\node[above, align=center]
at (axis cs:3.8,54) {\small 54};

\end{axis}

\begin{axis}[%
set layers,
log origin=infty,
clip = false,
axis y line*=right,
every outer y axis line/.append style={mycolor2},
every y tick label/.append style={font=\footnotesize\color{mycolor2}},
every y tick/.append style={mycolor2},
ymode=log,
ymin=1e-8,
ymax=1,
ytick={1e-8, 1e-7, 1e-6, 1e-5, 1e-4, 1e-3, 1e-2, 1e-1, 1},
yminorticks=true,
ylabel style={font=\footnotesize\color{mycolor2}},
ylabel={global sensitivity indices $\varrho(\J_{\ell},\hat{u}_{\boldx_g})$},
yticklabel pos=right,
ymajorgrids,
]

\addplot[ybar, bar width=0.35, fill=mycolor2, draw=black, area legend, postaction={line space=12pt, pattern=my north east lines}, forget plot] table[row sep=crcr] {%
1.2	0.998921070745274\\
2.2	0.0521691082873199\\
3.2	0.00117146455316547\\
4.2	1.72520509733178e-05\\
};
\node[above, align=center]
at (axis cs:1.2,1) {\tiny 0.9989};
\node[above, align=center]
at (axis cs:2.2,0.0522) {\tiny 0.0522};
\node[above, align=center]
at (axis cs:3.2,0.0012) {\tiny 0.0012};
\node[above, align=center]
at (axis cs:4.2,1.7e-5) {\tiny 1.73e-05};

\addplot[ybar, bar width=0.35, fill=mycolor3, draw=black, area legend, forget plot] table[row sep=crcr] {%
1.2	0.947532645868731\\
2.2	0.00107374984318008\\
3.2	3.27179856117439e-06\\
4.2	4.2809053974168e-08\\
};
\node[below, align=center]
at (axis cs:1.2,0.9475) {\tiny 0.9475};
\node[below, align=center]
at (axis cs:2.2,0.0011) {\tiny 0.0011};
\node[below, align=center]
at (axis cs:3.2,3.2e-06) {\tiny 3.27e-06};
\node[below, align=center]
at (axis cs:4.2,4.2e-08) {\tiny 4.28e-08};

\addplot [color=black, dashed, forget plot]
  table[row sep=crcr]{%
1.025	0.985796282046403\\
1.375	0.985796282046403\\
};
\addplot [color=black, dotted, forget plot]
  table[row sep=crcr]{%
1.025	0.990217211513946\\
1.375	0.990217211513946\\
};

\addplot [color=black, dashed, line width=0.5mm]
  table[row sep=crcr]{%
2.025	0.0141048316491327\\
2.375	0.0141048316491327\\
};
\addplot [color=black, dotted, line width=0.5mm]
  table[row sep=crcr]{%
2.025	0.00970928785122194\\
2.375	0.00970928785122194\\
};
\legend{mean, median}
\addplot [color=black, dashed, line width=0.5mm, forget plot]
  table[row sep=crcr]{%
3.025	9.8370828718121e-05\\
3.375	9.8370828718121e-05\\
};
\addplot [color=black, dotted, line width=0.5mm, forget plot]
  table[row sep=crcr]{%
3.025	4.31553234533777e-05\\
3.375	4.31553234533777e-05\\
};
\addplot [color=black, dashed, line width=0.5mm, forget plot]
  table[row sep=crcr]{%
4.025	5.15475744786508e-07\\
4.375	5.15475744786508e-07\\
};
\addplot [color=black, dotted, line width=0.5mm, forget plot]
  table[row sep=crcr]{%
4.025	1.27416663127499e-07\\
4.375	1.27416663127499e-07\\
};

\end{axis}
\end{tikzpicture}%
\caption{Cardinality (left, blue, solid) of the index sets $\J_\ell$ and the corresponding largest (right, orange, striped), smallest (right, yellow, solid), mean (dashed line) and median (dotted line) of the global sensitivity indices $\varrho(\J_\ell,u_{\boldx_g}^{\texttt{usFFT}})$ w.r.t.\ $\boldx_g$ for the affine example with $\eta=\textnormal{IX}$, i.e., $s=2000,\; N=32$.}
  \label{fig:affine_anova}
\end{figure}
Figure \ref{fig:affine_anova} again illustrates the cardinality of the sets $\J_\ell$ as well as their GSI.

\subsubsection*{Discussion}

We note that even for small sparsity parameters $s$ the approximation $\E(u_{\boldx_g}^{\texttt{usFFT}})$ seems to be quite accurate in Figure \ref{fig:affine_ew}. Unfortunately, for all other parameter choices $\eta$ except I and II, the  magnitude of the errors does not decrease any further than $1.5 \cdot 10^{-5}$, which is probably caused by the poor performance of the Monte-Carlo approximation $\overline{\check{u}_{\boldx_g}}$.

The magnitudes of the errors $\err_2^{\eta}$ in Figure \ref{fig:affine_decay} are already very low for small sparsity parameters $s$ compared to the expected function values shown in Figure \ref{fig:affine_ew1}. We note that there is an obvious improvement for each sparsity parameter $s$, when we progress from $N=32$, the first data point in each cluster, to $N=64$, the second data point. In the periodic case, the important frequencies are very well localized around zero, such that the choice of $N$ had almost no impact. This time, we really lose some accuracy if we choose the smaller extension $N=32$. For the sparsity parameter $s=2000$ we also see this effect when progressing from $N=64$ to $N=128$. The overall decay of the error is a lot slower compared to the periodic example. This is probably mainly caused by the non-smooth tent transformation used, cf.\ Section \ref{subsubsec:affine}. Again, some numerical tests determining the error $\text{err}_\infty^{\eta}$ revealed a similar behavior as in the periodic setting and showed that these errors again are not larger than at most $20$ times the error $\err_2^{\eta}$.

Since the Fourier coefficients do not decay as fast as in the smooth periodic case, we detected a significantly larger number of one- and two-dimensional couplings in Figure \ref{fig:affine_anova}. Again, the frequencies with only one non-zero entry explain the largest part of the variance of the function, but this time the minimum percentage is lower than in the periodic example with only about $94.5\%$. Accordingly, the importance of the two- and three-dimensional pairings did slightly grow. The large number of important coefficients with only one, two or three non-zero entries also results in nearly no detected significant frequencies with more than three non-zero entries. For example, the $54$ frequencies in $\J_4$ will vanish when working with larger extensions $N$, since other frequencies with less entries are preferred in that case. So even though we are working with the moderate stochastic dimension $d = 20$, we do not detect any frequencies, where the half or even only a quarter of these dimensions are active simultaneously.

\subsection{Lognormal example} \label{subsec:lognormal}

We consider a two-dimensional problem based on the example in \cite{CheHoYaZha13} on the domain $D_{\boldx}=[0,1]^2$ with right-hand side $f(\boldx) = \sin(1.3\pi x_1 + 3.4 \pi x_2) \cos(4.3 \pi x_1 - 3.1 \pi x_2)$. The lognormal random coefficient is given by
\begin{align*}
a(\boldx,\boldy) \coloneqq \exp(b(\boldx,\boldy)) \qquad \text{and} \qquad b(\boldx,\boldy) \coloneqq \sum_{j=1}^{d} \frac{1}{j} y_j \psi_j(\boldx)
\end{align*}
with the functions
\begin{align*}
\psi_j(\boldx) \coloneqq \sin(2\pi j x_1)\cos(2 \pi (d+1 - j) x_2).
\end{align*}
In \cite{CheHoYaZha13}, the stochastic dimension $d=4$ has been used. Here, we will work with $d = 10$ to receive a more complicated and higher-dimensional problem setting. We use a standard normally distributed random variable $\boldy \sim \mathcal{N}(\boldzero,\boldI)$ with $\boldI$ the identity matrix of dimension $d$ as before. Hence, we have, that for each $\boldx$ there holds $0 < a(\boldx,\boldy) < \infty$ for any $\boldy$. However, there do not exist the constants $0 < a_{\min} \leq a_{\max}<\infty$ in this example, since $b(\boldx,\boldy)$ can become arbitrarily small or large. Therefore, the problem is neither uniformly elliptic nor uniformly bounded. This complicates the analysis of this problem tremendously. We can still stick with it for our numerical tests, since we only need the solvability of the differential equation for fixed values of $\boldy$. Further, we have $b(\boldx,\boldy) \in [-3,3]$ and therefore $\exp(b(\boldx,\boldy)) \in [\e^{-3},\e^3] \approx [0.05, 20.09]$ with a probability of more than $99\%$ for each $\boldx \in D_{\boldx}$, i.e., tremendously small or large values of $a(\boldx,\boldy)$ are very unlikely to appear.

\begin{figure}[tb]
	\centering	
	\subfloat[$\overline{\check{u}_{\boldx_g}}$]{
		\includegraphics[scale=0.35]{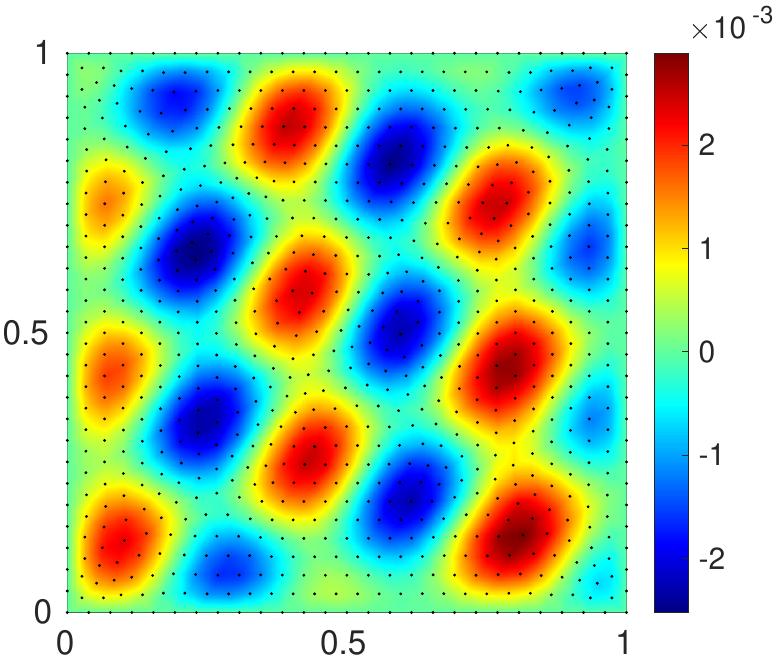} \label{fig:lognormal_ew1}
	}
	~	
	\subfloat[$\eta =$ IV ($s=500, N=32$)]{
    		\includegraphics[scale=0.35]{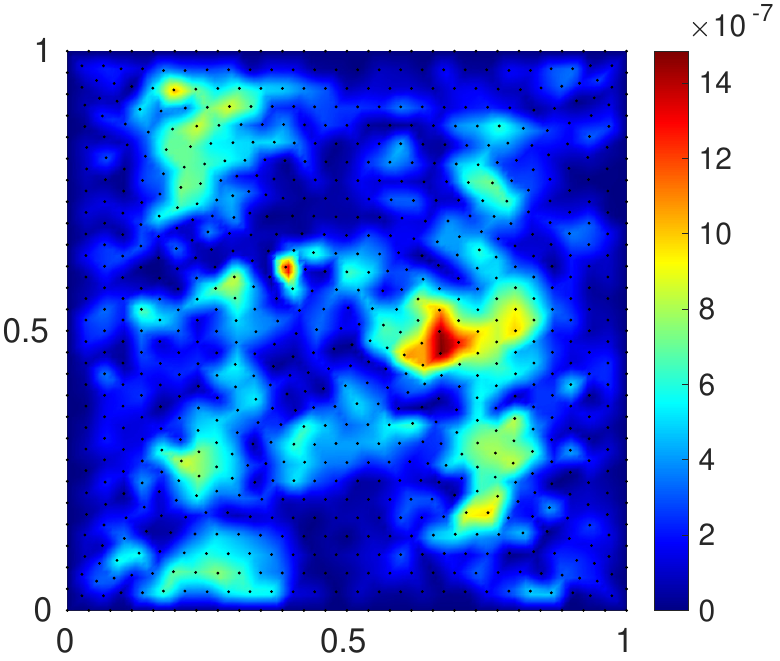} 
	}
	~
	\subfloat[$\eta =$ VI ($s=1000, N=32$)]{
    		\includegraphics[scale=0.35]{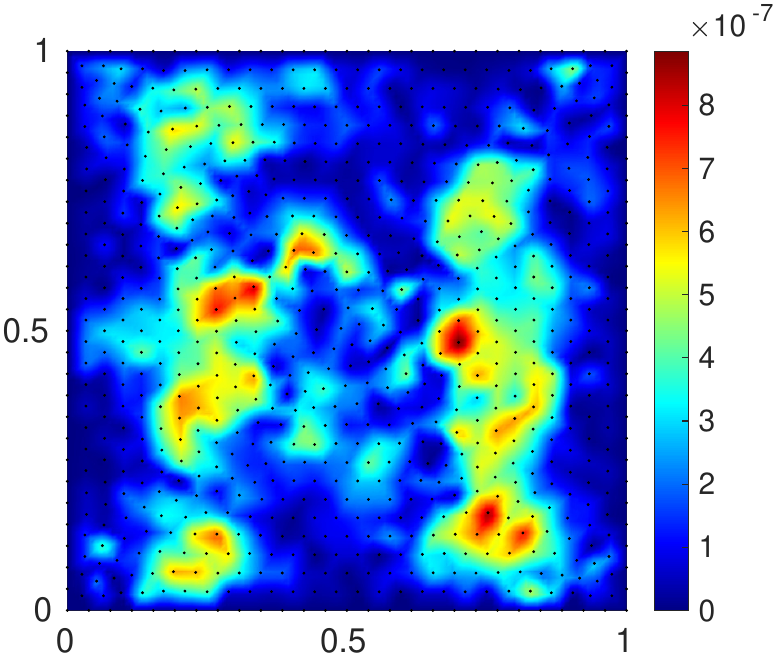} 
	}
	\caption{The MC approximation $\overline{\check{u}_{\boldx_g}}$ and the pointwise errors $\vert \overline{\check{u}_{\boldx_g}} - \E(u_{\boldx_g}^{\texttt{usFFT}})\vert$ for $\eta =$ IV and VI for the lognormal example.}
  \label{fig:lognormal_ew}
\end{figure}

Figure \ref{fig:lognormal_ew1} once again illustrates the  Monte-Carlo approximation of the expectation value $\overline{\check{u}_{\boldx_g}}$ with $n_{\text{MC}} = 10^6$ samples used. 
\begin{figure}[tb]
	\centering	
    	\includegraphics[scale=0.65]{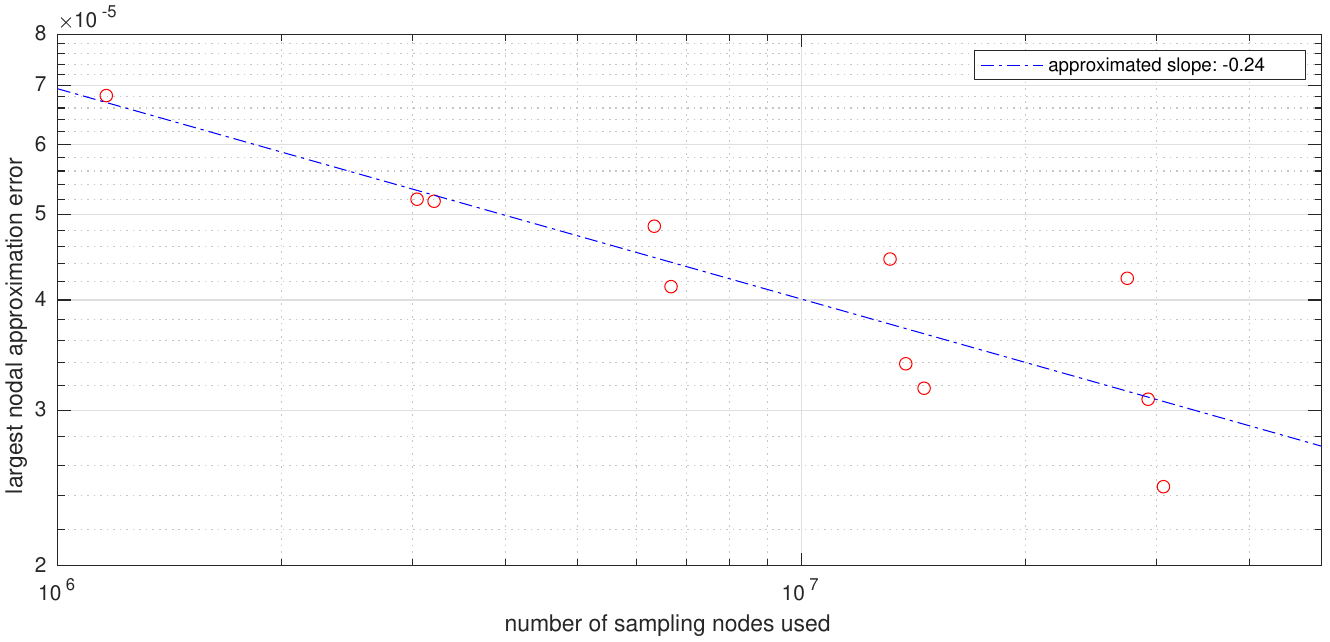} 
	\caption{Largest error $\err_2^{\eta}$ w.r.t.\ the nodes $\boldx_g$ for the parameter settings $\eta=$ I to XI displayed in Table \ref{table:parameters} for the lognormal example.}
  \label{fig:lognormal_decay}
\end{figure}
The decay of the largest error $\err_2^{\eta}$ w.r.t.\ the nodes $\boldx_g$ is shown in Figure \ref{fig:lognormal_decay}, where the data points are ordered from left to right w.r.t.\ increasing $\eta$ as in Figure \ref{fig:affine_decay}.
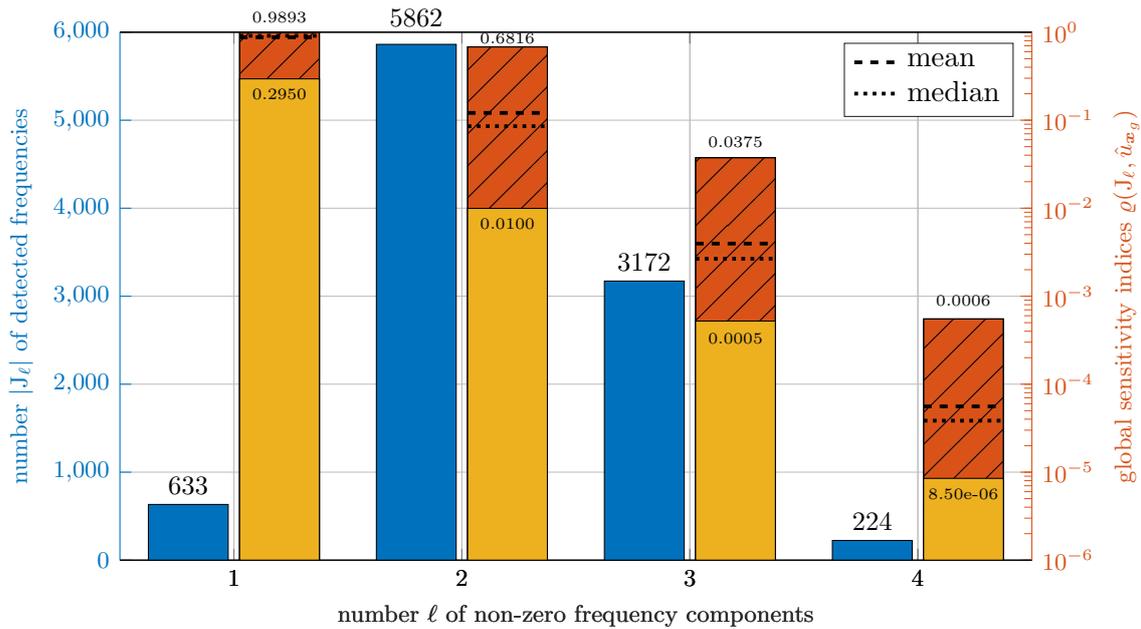
\begin{figure}[tb]
	\centering	
	\setlength\fwidth{12cm}
	\setlength\fheight{7cm}
\definecolor{mycolor1}{rgb}{0.00000,0.44700,0.74100}%
\definecolor{mycolor2}{rgb}{0.85000,0.32500,0.09800}%
\definecolor{mycolor3}{rgb}{0.92900,0.69400,0.12500}%
\begin{tikzpicture}
\pgfplotsset{
	scale only axis,
	xmin = 0.5,
	xmax = 4.5,
	xtick={1, 2, 3, 4},
	every x tick label/.append style={font=\footnotesize},
	xlabel style={font=\footnotesize\color{white!15!black}},
	xlabel={number $\ell$ of non-zero frequency components},
	xmajorgrids,
	width=\fwidth,
	height=\fheight,
	at={(0\fwidth,0\fheight)},
	legend style={legend cell align=left, align=left, draw=white!15!black}
}

\begin{axis}[
set layers,
clip = false,
axis y line*=left,
every outer y axis line/.append style={mycolor1},
every y tick label/.append style={font=\footnotesize\color{mycolor1}},
every y tick/.append style={mycolor1},
ymode=linear,
ymin=0,
ymax=6000,
ylabel style={font=\footnotesize\color{mycolor1}},
ylabel={number $\vert\J_{\ell}\vert$ of detected frequencies},
yticklabel pos=left,
]

\addplot[ybar, bar width=0.35, fill=mycolor1, draw=black, area legend] table[row sep=crcr] {%
0.8	633\\
1.8	5862\\
2.8	3172\\
3.8	224\\
};
\node[above, align=center]
at (axis cs:0.8,633) {\small 633};
\node[above, align=center]
at (axis cs:1.8,5950) {\small 5862};
\node[above, align=center]
at (axis cs:2.8,3172) {\small 3172};
\node[above, align=center]
at (axis cs:3.8,224) {\small 224};

\end{axis}

\begin{axis}[%
set layers,
log origin=infty,
clip = false,
axis y line*=right,
every outer y axis line/.append style={mycolor2},
every y tick label/.append style={font=\footnotesize\color{mycolor2}},
every y tick/.append style={mycolor2},
ymode=log,
ymin=1e-6,
ymax=1,
ytick={1e-6, 1e-5, 1e-4, 1e-3, 1e-2, 1e-1, 1},
yminorticks=true,
ylabel style={font=\footnotesize\color{mycolor2}},
ylabel={global sensitivity indices $\varrho(\J_{\ell},\hat{u}_{\boldx_g})$},
yticklabel pos=right,
ymajorgrids,
]

\addplot[ybar, bar width=0.35, fill=mycolor2, draw=black, area legend, postaction={line space=12pt, pattern=my north east lines}, forget plot] table[row sep=crcr] {%
1.2	0.989300954163508\\
2.2	0.681564887018414\\
3.2	0.0375033078739069\\
4.2	0.0005517079158456\\
};
\node[above, align=center]
at (axis cs:1.2,0.9893) {\tiny 0.9893};
\node[above, align=center]
at (axis cs:2.2,0.58) {\tiny 0.6816};
\node[above, align=center]
at (axis cs:3.2,0.0375) {\tiny 0.0375};
\node[above, align=center]
at (axis cs:4.2,0.0006) {\tiny 0.0006};

\addplot[ybar, bar width=0.35, fill=mycolor3, draw=black, area legend, forget plot] table[row sep=crcr] {%
1.2	0.295030570940103\\
2.2	0.00995188892231512\\
3.2	0.000521463742375137\\
4.2	8.50187391125605e-06\\
};
\node[below, align=center]
at (axis cs:1.2,0.2950) {\tiny 0.2950};
\node[below, align=center]
at (axis cs:2.2,0.0100) {\tiny 0.0100};
\node[below, align=center]
at (axis cs:3.2,0.0005) {\tiny 0.0005};
\node[below, align=center]
at (axis cs:4.2,8.50e-06) {\tiny 8.50e-06};

\addplot [color=black, dashed, line width=0.5mm, forget plot]
  table[row sep=crcr]{%
1.025	0.874959015726121\\
1.375	0.874959015726121\\
};
\addplot [color=black, dotted, line width=0.5mm, forget plot]
  table[row sep=crcr]{%
1.025	0.911216864511522\\
1.375	0.911216864511522\\
};

\addplot [color=black, dashed, line width=0.5mm]
  table[row sep=crcr]{%
2.025	0.121025095306034\\
2.375	0.121025095306034\\
};
\addplot [color=black, dotted, line width=0.5mm]
  table[row sep=crcr]{%
2.025	0.085442375588141\\
2.375	0.085442375588141\\
};
\legend{mean, median}
\addplot [color=black, dashed, line width=0.5mm, forget plot]
  table[row sep=crcr]{%
3.025	0.00395999351651479\\
3.375	0.00395999351651479\\
};
\addplot [color=black, dotted, line width=0.5mm, forget plot]
  table[row sep=crcr]{%
3.025	0.00266166020221721\\
3.375	0.00266166020221721\\
};
\addplot [color=black, dashed, line width=0.5mm, forget plot]
  table[row sep=crcr]{%
4.025	5.58954513296413e-05\\
4.375	5.58954513296413e-05\\
};
\addplot [color=black, dotted, line width=0.5mm, forget plot]
  table[row sep=crcr]{%
4.025	3.84855949459161e-05\\
4.375	3.84855949459161e-05\\
};

\end{axis}
\end{tikzpicture}%
\caption{Cardinality (left, blue, solid) of the index sets $\J_\ell$ and the corresponding largest (right, orange, striped), smallest (right, yellow, solid), mean (dashed line) and median (dotted line) of the global sensitivity indices $\varrho(\J_\ell,u_{\boldx_g}^{\texttt{usFFT}})$ w.r.t.\ $\boldx_g$ for the lognormal example with $s=2000,\; N=32$.}
  \label{fig:lognormal_anova}
\end{figure}
Finally, Figure \ref{fig:lognormal_anova} shows the cardinality of the sets $\J_\ell$ as well as their GSI for this example.

\subsubsection*{Discussion}

We note that the pointwise solution in Figure \ref{fig:lognormal_ew1} has a more interesting structure than for the other examples above, mainly caused by the lognormal random coefficient and the non-constant right-hand side $f(\boldx)$. Nevertheless, the approximations $\E(u_{\boldx_g}^{\texttt{usFFT}})$ achieve small errors, which are shown in Figure \ref{fig:lognormal_ew}. This time, a further increase of the sparsity parameter $s$ and the extension $N$ still increase the accuracy of our approximations, so the stagnation due to the limitations of the Monte-Carlo approximation $\overline{\check{u}_{\boldx_g}}$, that we saw in the previous examples, does not occur yet.

The pointwise errors $\err_2^{\eta}(\boldx_g)$ behave slightly worse but still very good, as we see in Figure \ref{fig:lognormal_decay}. Again, the increase of the extension $N$ shows visible improvements of the approximation error $\err_2^{\eta}$. The decay rate is lower than before, matching our expectations since the lognormal example is far more difficult than the affine or periodic examples. Note that once again the slope considers all data points shown, while specific decays for fixed extensions $N$ might be slower or faster. Further, the size of the error $\text{err}_\infty^{\eta}$ is again about $10$ times the size of $\err_2^{\eta}$, revealing also a good pointwise approximation w.r.t.\ the random variable $\boldy$ in this scenario.
                
We notice a similar distribution of the detected frequencies $\boldk$ to the index sets $\J_{\ell}\cap I$ as before, cf.\ Figure \ref{fig:lognormal_anova}. The key difference is the size of the GSI for each of these index sets. The range of the GSI for $\J_1$ increased significantly, the minimal portion of variance is now only about $30\%$. Obviously, the GSI for the other index sets $\J_\ell$ grew accordingly. This is probably caused by the more difficult structure of the lognormal diffusion coefficient $a$ and the corresponding more difficult structure of the solution which is reflected in larger differences in the optimal frequency sets $\I_{\boldx_g}$, $g=1,\ldots,G$, cf.\ Remark \ref{rem:growth}. Nevertheless, we again detect nearly no significant frequencies $\boldk$ with $4$ or more active dimensions as in the previous examples.
                
\begin{remark} \label{rem:growth}
As mentioned before, the output of the usFFT contains more than the sparsity parameter $s$ frequencies since we join the detected index sets in each dimension increment and use no thresholding technique to reduce the number of found frequencies after that. While we have no reasonably tight theoretical bounds on the size of the output yet, we can further investigate the number of output frequencies in our numerical tests. In detail, we express the detected output sparsity $s_{\text{real}}$ as a multiple of the given sparsity parameter $s$, i.e., $s_{\text{real}} = q \cdot s$ with some factor $q \in \R$.

In the numerical tests for the first periodic example in Section \ref{subsec:periodic}, i.e., $\mu=1.2$ and $c=0.4$, we have $q \in [2.41, 2.74]$, where the larger values of $q$ tend to appear for smaller sparsity parameters $s$. For the quickly decaying example, i.e., $\mu=3.6$ and $c=1.5$, we have $q \in [1.9042, 2.45]$ and again the larger values of $q$ are attained for small sparsity parameters $s$.

The affine model in Section \ref{subsec:affine} results in $q \in [2.06, 2.186]$, where $q < 2.1$ is only attained for $\eta =$ I, II and IV, so parameter settings with small sparsity parameters $s$ and extension $N=32$.

Finally, in the complicated lognormal case in Section \ref{subsec:lognormal}, we observe $q \in [4.776, 5.15]$. While the values above $5$ only appear for $\eta =$ I and II, we still have significantly larger factors $q$ than before%
.
However, the magnitude of $q$ is still very small compared to the size $\vert\setT_G\vert = G = 739$ in our examples.

Our observation is consistent with recent results presented in \cite{KaKaKu20} which considers the periodic model only. The crucial common feature is that the pointwise approximations $u_{\boldx_g}$  can be regarded as elements of a joint reproducing kernel Hilbert space with uniformly good kernel approximants.

Overall, the factor $q$ in our examples is much smaller than $G$, which would be the worst factor possible in the case that all $\I_{\boldx_g}$ are disjoint. Hence, as already mentioned in Section \ref{sec:intro}, the given complexities in Theorem \ref{thm:complexities} are way too pessimistic and the true amount of sampling locations and computational steps needed is much smaller in all of our numerical examples.
\end{remark}

\subsection{Comparison to given frequency sets}
\label{subsec:known_freq}

The main effort of the usFFT lies in detecting the index set $\I \subset \Gamma$. The computation of the corresponding Fourier coefficients in the final step of Algorithm \ref{algo:sfft_general} needs significantly less samples than the detection steps before. Hence, the question arises, if an a priori choice of the index set $\I$ should be preferred to reduce the computational cost, cf.\ Remark \ref{rem:sizes}. Therefore, we now consider the following kinds of index sets:

\begin{itemize}
	\item	axis cross with uniform weight $1$: $ \I = \lbrace \boldk \in \Z^{d}: \|\boldk\|_0 = 1, \|\boldk\|_1 \leq N \rbrace$
	\item	hyperbolic cross with uniform weight $\frac14$: $ \I = \lbrace \boldk \in \Z^{d}: \prod_{j=1}^{d} \max(1,4 \vert k_j \vert) \leq N \rbrace$
	\item	hyperbolic cross with slowly or quickly $(q=1$ or $2)$ decaying weights $\frac{1}{j^q}$: $ \I = \lbrace \boldk \in \Z^{d}: \prod_{j=1}^{d} \max(1,j^q \vert k_j \vert) \leq N \rbrace$
	\item	$l_1$-ball with slowly or quickly $(q=1$ or $2)$ decaying weights $\frac{1}{j^q}$: $ \I = \lbrace \boldk \in \Z^{d}: \sum_{j=1}^{d} j^q \vert k_j \vert \leq N \rbrace$
\end{itemize}

\begin{figure}[tb]
	\centering	
	\subfloat[periodic example with $\mu = 1.2, c=0.4$]{
		\includegraphics[scale=0.58]{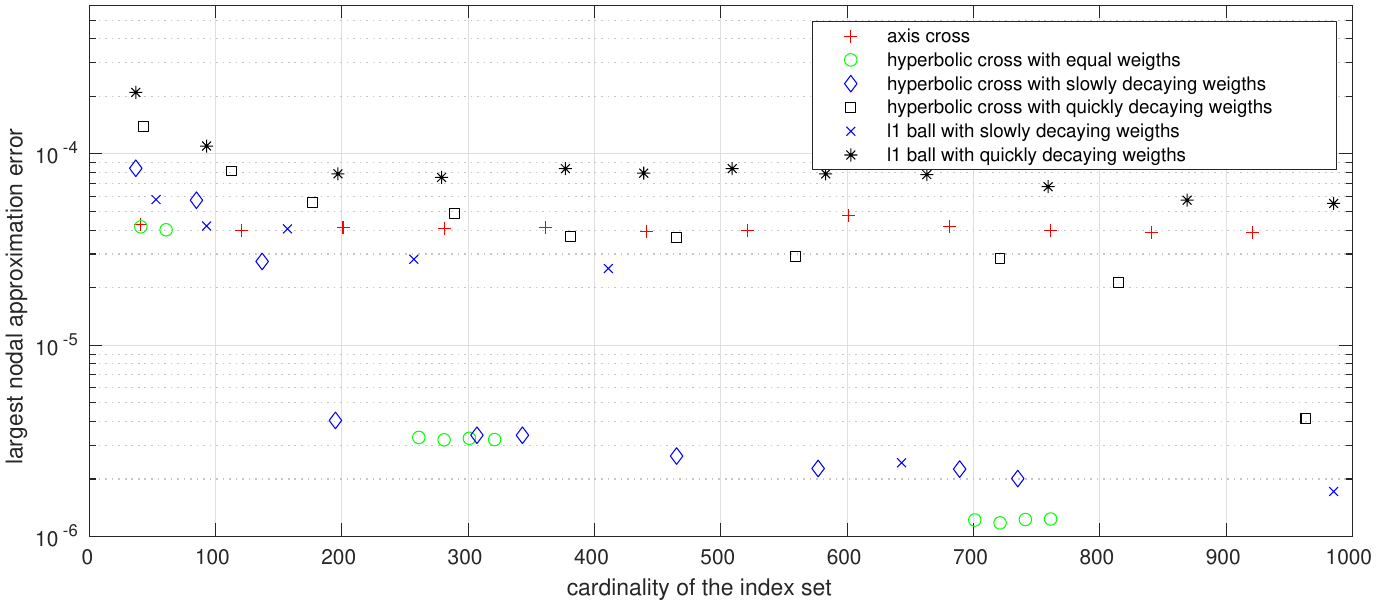} 
	}
	\\	
	\subfloat[periodic example with $\mu = 3.6, c=1.5$]{
		\includegraphics[scale=0.58]{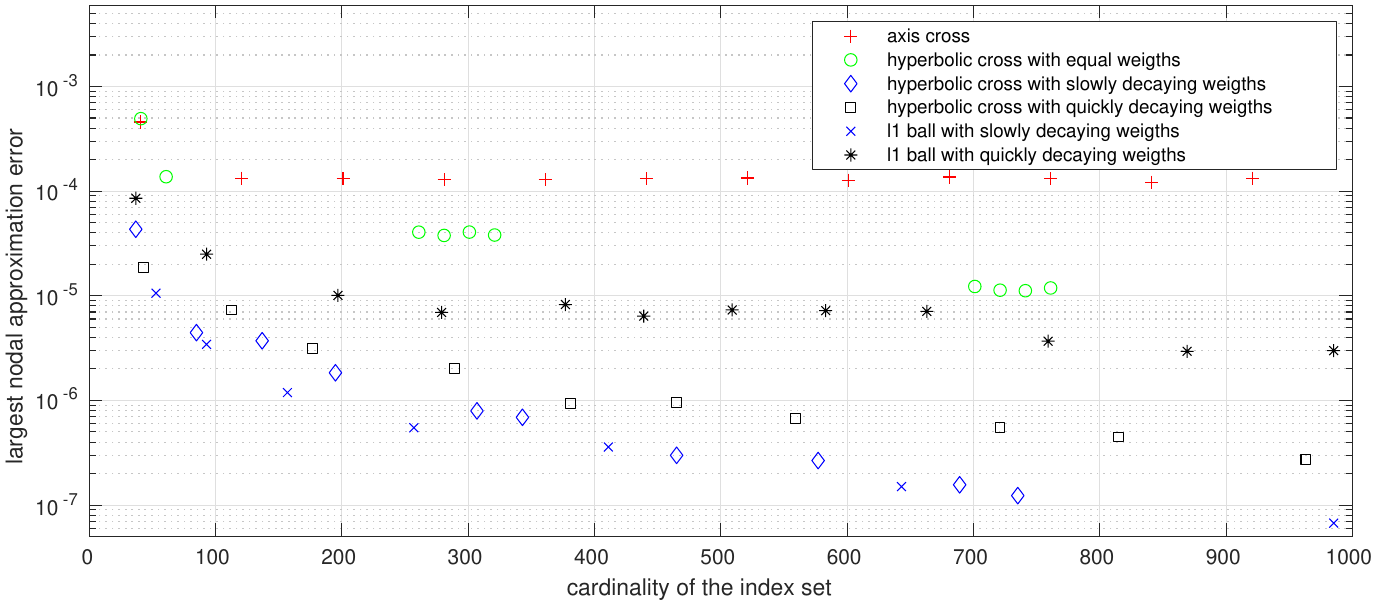} 
	}
	\\	
	\subfloat[affine example]{
    		\includegraphics[scale=0.58]{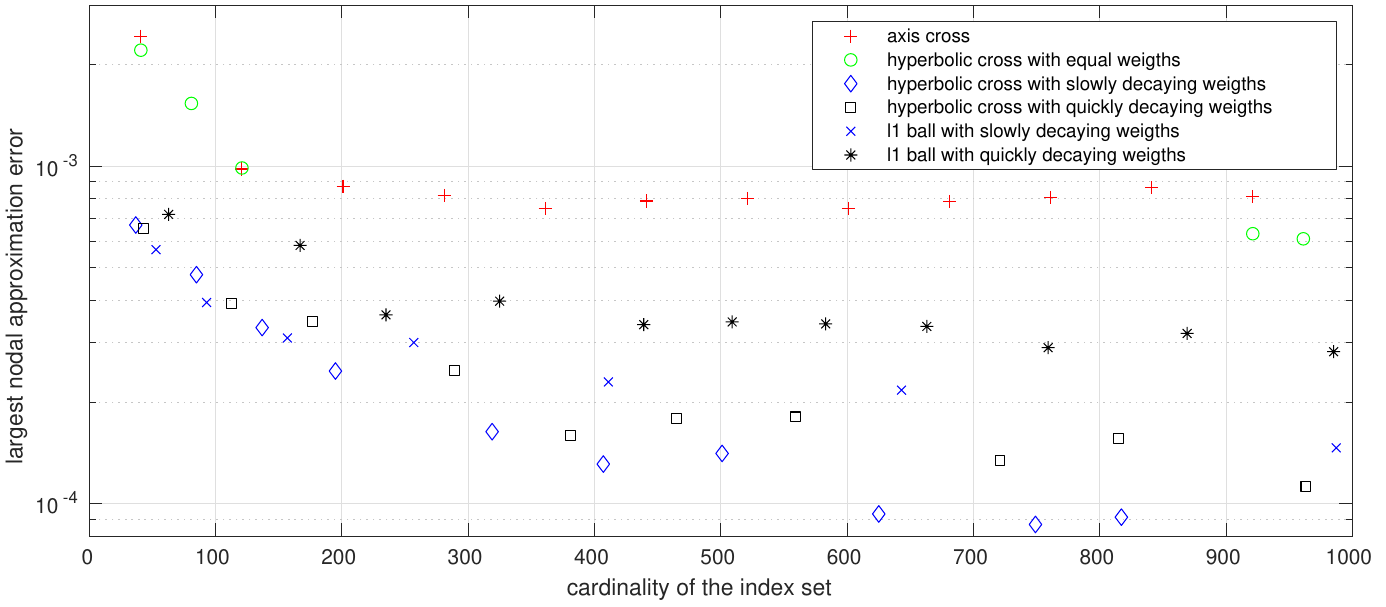} 
	}
	\caption{Largest error $\err_2^{\eta}$ w.r.t.\ the nodes $\boldx_g$ for the periodic and affine examples with given frequency sets.}
  \label{fig:freq_known}
\end{figure}

The Fourier coefficients $c_{\boldk}(u_{\boldx_g})$ are approximated using the same multiple R1L approach as in step~3 of our Algorithm \ref{algo:sfft_general}, i.e., we just skipped steps~1 and~2 by choosing the index set $\I$ instead of detecting it. Figure \ref{fig:freq_known} illustrates the largest error $\err_2^{\eta}$ w.r.t.\ the nodes $\boldx_g\in\mathcal{T}_G$ for the previously considered periodic and affine examples, cf.\ Sections~\ref{subsec:periodic} and~\ref{subsec:affine}, with these given frequency sets $\I$ for various refinements $N$. 

The magnitude of the errors is considerably larger than for comparable parameter settings of the usFFT, e.g., $\eta =$ I to III, especially for the periodic example. Further, we also see that the particular choice of the structure of the index set plays an important role. Obviously, a cleverly chosen index set reduces the size of the approximation error tremendously, especially in the periodic settings. But finding a good or even optimal choice of the index set is highly non-trivial, since it requires sufficient a priori information about the PDE and the structure of its solution or additional computational effort, e.g., to determine suitable weights for a given index set structure. This can be observed for example when comparing the hyperbolic cross index sets for the periodic examples. The uniform weights achieve the best results for $\mu = 1.2$, but cannot keep up at all with the decaying weights for the faster decay rate $\mu = 3.6$. On the other hand, even if we know, that there is a certain decay in our random coefficient, it is not clear how to choose suitable decay rates for the weights in order to guarantee reasonable results -- specifically in pre-asymptotic settings, which is the rule rather than the exception when numerically determining solutions of high-dimensional problems.

The usFFT does not depend on these kind of information, as its choice of the frequency set is fully adaptive and the only required a priori information is the search space $\Gamma$, which can be chosen sufficiently large without disturbing the results of the algorithm. Further, the detected frequency set $\I$ provides these additional information about the structure of the solution $u$ as well as the dependence on the random variables $\boldy$. In other words, the additional amount of samples needed for the usFFT makes these structural information unnecessary, detects them on its own and provides a possibility to extract them afterwards from the output.

\begin{remark}\label{rem:sizes}
The computations in this section and in step~3 of the usFFT are performed using the multiple R1L approach for the efficient computation of Fourier coefficients for a given frequency set $\I$ as proposed in \cite{Kae17}. From \cite[Cor.~3.7]{Kae17} we get a bound on the number of sampling nodes $M$ used. Since we are working with $c=2$ and $\delta = 0.5$ as in \cite[Alg.~3]{KaPoVo17}, we arrive at $M \leq \lceil 2 \,\ln (2\vert \I \vert) \rceil 4 (\vert \I \vert -1)$. Note that this upper bound is very rough and the actual number of used sampling nodes in almost all numerical experiments is much lower.

As stated above several times, this number of samples used in step~3 of the usFFT is just a small fraction of the total number of used samples when applying the usFFT. In particular, the computation of the actual Fourier coefficients $c_{\boldk}^{\textnormal{\texttt{usFFT}}}(u_{\boldx_g})$ for the detected frequency set $\I$ requires roughly $0.4\%$ or $0.3\%$ of the total sampling amount for the two different parameter choices of the periodic example in Section \ref{subsec:periodic}, around $0.1\%$ in the affine case in Section \ref{subsec:affine} and about $0.65\%$ for the lognormal model from Section \ref{subsec:lognormal}.
\end{remark}

In \cite{CheHoYaZha13,ZhaHuHoLiYa14}, a data-driven method was proposed, which is capable of computing approximations for multiple right-hand sides $f(\boldx)$ from a certain function class. Our usFFT approach can also be generalized in a similar, data-driven way: For a given class of functions $f(\boldx)$ or even $f(\boldx,\boldy)$, we can use the usFFT in order to compute the frequency set $\I$ for one randomly selected right-hand side $f$ or randomly select multiple right-hand sides $f$ and compute unions of the corresponding index sets $\I_f$ by means of (a slight modification of) the presented usFFT. In each case, we end up with a frequency set $\I$, which is probably a good choice for all the functions $f$ in the given class, since they are hopefully very similar to each other. Hence, we can use this index set $\I$ as a starting point and compute approximations of the corresponding Fourier coefficients $c_{\boldk}(u_{\boldx_g}), \boldk \in \I,$ as done above. This approximation of $u_{\boldx_g}$ is then probably a lot better, i.e., the detected index set $\I$ is a better localization of the largest Fourier coefficients than some a priori choice.

\section{Summary}
The proposed dimension-incremental method provides an efficient possibility to adaptively compute solutions to parametric PDEs involving high-dimensional random coefficients. The non-intrusive behavior allows the use of different, suitable PDE solvers and therefore generates high adaptability of our algorithm to different problem settings as well. We show, that the amount of PDE solutions needed can be decreased by our method compared to the naive repetitive approach even in the worst case. The numerical experiments underline this and show the functionality of our method for practical examples.

\section*{Acknowledegement}
L. K\"ammerer gratefully acknowledges funding by the Deutsche \linebreak Forschungsgemeinschaft (DFG, German Research Foundation) with the project
number \linebreak 380648269 and Daniel Potts with the project number 416228727 -- SFB 1410.

\begin{appendices}

\section{Rank-1 lattices in the sFFT}\label{app:r1l}
Popular approaches with the properties stated at the end of Section \ref{subsec:method} are for example based on so-called rank-1 lattices. A rank-1 lattice (R1L) is a set 
\begin{align*}
\Lambda(\boldz,M) \coloneqq \left\lbrace \frac{i}{M} \boldz \text{ mod } \boldone:\, i=0,...,M-1 \right\rbrace
\end{align*}
with a so-called generating vector $\boldz \in \Z^d$ and lattice size $M \in \N$. In \cite{PoVo14}, single rank-1 lattices (single R1Ls) were used as sampling strategy in the dimension-incremental method and provided a perfectly stable, reliable and efficient way to reconstruct the projected Fourier coefficients $\tilde{\hat{p}}_{(1,...,t),\boldk}$. In \cite{KaPoVo17} and \cite{KaKrVo20}, other approaches based on multiple rank-1 lattices (multiple R1Ls) and random rank-1 lattices (random R1Ls) have been studied. The main advantage of these approaches is a smaller size and a significantly faster construction of the involved sampling sets $\X$, but therefore they involve some failure probability, which is not needed for the perfectly stable single R1L approach. Table~\ref{table:complexities} shows the sampling and arithmetic complexities of the dimension-incremental method when using these different sampling strategies based on R1Ls. Further notes on these approaches and their behavior when used in the dimension-incremental method can be found in the referred works.

\begin{table}[tb]
\centering
\begin{tabular}{c|c|c} 
 & sample complexity & computational complexity \\ 
 \hline
 single R1Ls\;\cite{PoVo14} & $\mathcal{O}(d r^3 s^2 N)$ & $\mathcal{O}(d r^3 s^3 + d r^3 s^2 N \, \log^{\mathcal{O}(1)}(...)) $ \\
 \hline
 multiple R1Ls\;\cite{KaPoVo17} & $\mathcal{O}(d r^2 s N \, \log^{\mathcal{O}(1)}(...))$ & $\mathcal{O}(d^2 r^2 s N \,\log^{\mathcal{O}(1)}(...))$ \\
 \hline
 random R1Ls\;\cite{KaKrVo20} & $\mathcal{O}(d r s \, \log^{\mathcal{O}(1)}(...)) $ & $\mathcal{O}(d^2 r s N \,\log^{\mathcal{O}(1)}(...)) $ \\
\end{tabular}
\caption{Sampling and arithmetic complexities of the sFFT approach (with high probability, cf.\ \cite{KaPoVo17, KaKrVo20}) when using different sampling strategies based on R1Ls, where $\Gamma\subset[-N,N]^d$, $s\ge|\operatorname{supp}\hat{p}|$, and $r$ is the number of random projections computed in each dimension-incremental step.}
\label{table:complexities}
\end{table}

\section{Proof of Theorem \ref{thm:complexities}}\label{app:proof}

\begin{proof}%
Note that in the following explanations as well as in the corresponding Tables \ref{table:proof1} and \ref{table:proof2} 'sample complexities' always refers to the cardinality of the set of sampling locations, since we assume the black box sampling algorithm to provide the samples for all $G$ trigonometric polynomials $p^{(g)}$ simultaneously. In addition, we assume $s_\mathrm{local}\lesssim s$.

Theorem~\ref{thm:complexities} is a slight modification of \cite[Thm.~2]{KaKrVo20}. To be more precise, we apply Algorithm~\ref{algo:sfft_general} using a random R1L approach in the role of Algorithm~$\operatorname{A}$ and a spatial discretization $\mathcal{Y}$ based on multiple R1Ls, cf.\ \cite{Kae17}, in step~3. Accordingly, we mainly refer to the analysis of the sample complexity and computational complexity of the sFFT using random R1Ls given in \cite[Sec.~3.2]{KaKrVo20} as well as the therefore necessary theoretical results from \cite[Sec.~4]{KaPoVo17}. 

The crucial difference to \cite[Thm.~2]{KaKrVo20} is that we have to take account of the modification that we demand for reconstructing not only one but even $G$ different trigonometric polynomials with possibly differing frequency supports $\I_g$, $g=1,\ldots,G$. In the following, we discuss the necessary modifications on the bounds and parameter choices discussed, proved and used in \cite[Sec.~3.2.2 and 3.2.3]{KaKrVo20} and \cite[Lem.~4.4 and Thm.~4.6]{KaPoVo17}, such that the corresponding results hold.  

The modifications when considering the usFFT can be separated in two different parts, first the possibly larger sets $\J_t$ and $\I^{(1,\ldots,d)}$ in steps~2 and~3, respectively, and second the modified failure probability. Note that $r, \gamma_{\operatorname{A}}$ and $\gamma_{\operatorname{B}}$ are the decisive parameters which provide estimates on the failure probability later on and, thus, both are discussed in the second part of the proof.

\subsubsection*{Part 1: Size of the frequency sets $\J_t$ and $\I^{(1,\ldots,d)}$}

We start with the candidate sets $\J_t$ in step~2b of Algorithm~\ref{algo:sfft_general} and observe, that the cardinality $\vert \J_t \vert$ of the set of frequency candidates 
$$
\J_t = (\I^{(1,...,t-1)} \times \I^{(t)}) \cap \P_{(1,...,t)}(\Gamma)\subset\left(\bigcup_{\substack{i=1,\ldots,\tilde{r}\\g=1,\ldots,G}}\tilde{J}_{t-1,i,g} \times \P_t(\Gamma)\right) \cap \P_{(1,...,t)}(\Gamma)
$$ in each dimension increment $t$ can be simply bounded by $|\J_t|\lesssim r\,s\,G\,N_\Gamma$, which contains an additional factor $G$ now.
Applying the sampling strategy suggested in \cite[Sec.~2.1]{KaPoVo17} together with \cite[Lem.~4.5]{KaPoVo17} directly yields $|\mathcal{X}_{t,i}|\lesssim\max(s,N_\Gamma)\log(|\J_t|/\gamma_{\operatorname{A}})\lesssim\max(s,N_\Gamma)\log(r\, s\, G\,N_\Gamma/\gamma_{\operatorname{A}})$ with $\gamma_{\operatorname{A}}$ the failure probability of Algorithm $\operatorname{A}$.

Further, the cardinality of the finally detected frequency set $\I^{(1,\ldots,d)}$ used in step~3 of our algorithm is bounded from above by $s\, G$ due to the same argumentation, since $\tilde{r}=1$ and $\tilde{s}=s$ when $t=d$ holds in step~2. Accordingly, we can apply \cite[Alg.~1]{KaPoVo17} in order to construct a spatial discretization $\mathcal{Y}$ of $\I^{(1,\ldots,d)}$ based on multiple R1Ls which has a cardinality bounded by $|\mathcal{Y}|\lesssim\max(s\,G,N_\Gamma)\log(s\,G/\gamma_{\operatorname{B}})$, where $\gamma_{\operatorname{B}}$ is the failure probability of the construction of this spatial discretization, cf.\ \cite[Thm.~4.1]{Kae17}.

\subsubsection*{Sample and computational complexity of the usFFT}

Now, we need to discuss the computational complexities of the individual steps of the usFFT. Obviously, step~1 applies $d\,r\,G$ different one-dimensional FFTs of lengths at most $N_\Gamma$, which yields a computational complexity in $\OO{d\,r\,G\,N_\Gamma\log(N_\Gamma)}$.
In step~2, we apply $((d-2)\,r+1)\,G$ times \cite[Alg.~4]{KaKrVo20}, where the sampling set is a union of $L\in\OO{\log(r\, s\,G\, N_\Gamma/\gamma_{\operatorname{A}})}$ R1Ls of size at most in $\OO{\max \{s,N_\Gamma \} }$ and the input set of frequencies $\J_t$ is bounded from above by $\OO{r\,s\,G\,N_\Gamma}$ in its cardinality, which yields an arithmetic complexity in 
\begin{equation*}
\begin{split}
&\OO{d\,r\,G(\max(s, N_{\Gamma})\log(s\, N_{\Gamma})+d\,r\,s\,G\,N_{\Gamma})\log(r\, s\,G\, N_\Gamma/\gamma_{\operatorname{A}})}\\
&\hspace*{20em}
\subset
\OO{d^2\,r^2\, s\,G^2\, N_{\Gamma}\log^2(r\, s\, N_\Gamma\, G/\gamma_{\operatorname{A}})}
\end{split}
\end{equation*}
in the worst case.
Moreover, we observe $|\J_t|\lesssim s\,G\,N_\Gamma$ with a certain probability since the signals $p$ are all trigonometric polynomials and in the case where Algorithm~$\operatorname{A}$ does not fail in any case, we have $\bigcup_{i=1,\ldots,\tilde{r}} \tilde{\J}_{t,i,g}\subset \I^{(1,\ldots,t)}_{g}$ with $|\I^{(1,\ldots,t)}_{g}|\le |\I_g^{(1,\ldots,d)}|\le s$. As a consequence, we save a linear $r$ and the $r$ in the log term compared to the worst case arithmetic complexity, cf.\ \cite[Sec.~3.2.2]{KaKrVo20} for a similar argumentation. In addition, the same argumentation saves a factor $r$ in the logarithmic term of the upper bound on the number of sampling locations in step~2 with the same probability. Later, we specifically choose the parameters $r, \gamma_{\operatorname{A}}$ and $\gamma_{\operatorname{B}}$ such that the estimates hold with high probability - for that reason, we call these complexities with high probability complexities already here.

For computing the $G$ FFTs of step~3 we apply \cite[Alg.~2]{KaPoVo17}, which yields a computational complexity in 
\begin{equation*}
\begin{split}
&\OO{G\log(s\,G/\gamma_{\operatorname{B}})\big(\max(s\,G,N_\Gamma)\log(s\,G\, N_\Gamma)+s\,G(d+\log(s\,G))\big)}\\
&\hspace*{15em}
\subset
\OO{G\,\max(s\,G,N_\Gamma)\,\log (s\, G/ \gamma_{\operatorname{B}})\, (d+\log(s\,G \,N_\Gamma))}.
\end{split}
\end{equation*}

The sample complexities and computational complexities of the usFFT due to these modifications are given in Table \ref{table:proof1}, see \cite[Tab.~3.2]{KaKrVo20} for comparison to the sFFT. Here, the only changes are several appearances of the parameter $G$, i.e., for $G=1$ we observe the complexities of the sFFT.

\begin{table}
\centering
\begin{tabular}{l l l}
& sample complexity & computational complexity \\
\hline
Step 1 & $d\,r\,N_\Gamma$ & $d\,r\,G\,N_\Gamma\,\log N_\Gamma$ \\
\rule{0em}{1.25em}Step 2 (w.h.p.) & $d\,r\,\max(s,N_\Gamma)\,\log \frac{s \,G \,N_\Gamma}{\gamma_A}$ & $d^2 \,r \,s \,G^2 \,N_\Gamma \,\log^2 \frac{s \,G \,N_\Gamma}{\gamma_A}$ \\
\rule{0em}{1.25em}Step 2 (w.c.) & $d\,r\,\max(s,N_\Gamma)\,\log \frac{r\, s \,G \,N_\Gamma}{\gamma_A}$ & $d^2 \,r^2 \,s \, G^2 \,N_\Gamma \,\log^2 \frac{r\,s \,G \,N_\Gamma}{\gamma_A}$ \\
\rule{0em}{1.25em}Step 3 & $\max(s\,G,N_\Gamma)\,\log \frac{s\, G}{\gamma_{\operatorname{B}}}$ & $G\,\max(s\,G,N_\Gamma)\,\log \frac{s\, G}{\gamma_{\operatorname{B}}}\, (d+\log(s\,G \,N_\Gamma))$
\end{tabular}
\caption{Sample complexities and computational complexities with high probability (w.h.p.) and in the worst case (w.c.) for the different steps of Algorithm \ref{algo:sfft_general}, where the efficient identification by \cite[Alg.~4]{KaKrVo20} is used in step~2 and the multiple R1L approach from \cite[Alg.~1]{KaPoVo17} in step~3.}\label{table:proof1}
\end{table}

\subsubsection*{Part 2: Parameter choices}

We continue with the aforementioned second big part, where we need to discuss suitable choices of $r$, $\gamma_{\operatorname{A}}$ and $\gamma_{\operatorname{B}}$ to obtain our desired failure probability $\delta$.

To this end, we first consider the projection failure probability, i.e., the failure that occurs if important projected Fourier coefficients are close to zero and, thus, not detectable. The number $r$ of detection iterations determines, how many of these projections are computed. The more projections are considered, the less is the probability that a specific projected Fourier coefficient is small for all of them and hence not detectable. Therefore, the parameter $r$ directly controls this projection failure probability.

\subsubsection*{The number $r$ of detection iterations}

We consider a single trigonometric polynomial $p\not\equiv 0$ with $\min_{\boldh\in\text{supp }{\hat{p}}}|\hat{p}_g|\ge 3\theta$ and $\Gamma\supset\text{supp }{\hat{p}}$, $|\text{supp }{\hat{p}}|\le s$. Choosing
\begin{align*}
r =\lceil 2s(\log 3 + \log d + \log s +\log G - \log \delta)\rceil,
\end{align*}
as given in \cite[Lem. 7]{KaPoVo17}, yields a probability of at most $\frac{\delta}{3\,d\,s\,G}$ that all the projected Fourier coefficients are less than $\theta$ for at least one frequency.

For $G$ different of such trigonometric polynomials $p^{(g)}$, we then apply the union bound. Therefore, the probability, that all the projected Fourier coefficients are less than $\theta$ for at least one frequency and at least one signal, is bounded by $\frac{\delta}{3\,d\,s}$.

\subsubsection*{The failure probabilities $\gamma_{\operatorname{A}}$ and $\gamma_{\operatorname{B}}$}

In the remaining lines of Part 2, we investigate the choices of the failure probabilities $\gamma_{\operatorname{A}}$ and $\gamma_{\operatorname{B}}$.

We start with the parameter $\gamma_{\operatorname{A}}$, which is in fact the failure probability of \cite[Alg.~4]{KaKrVo20} in the role of Algorithm~$\operatorname{A}$. When choosing $\gamma_{\operatorname{A}}:=\frac{\delta}{3\,d\,s\,G}$, we observe, that the probability, that at least one of the $G$ applications of Algorithm~$\operatorname{A}$ fails in step~2d (for fixed $t$ and $i$), is bounded from above due to the union bound by $\frac{\delta}{3\,d\,s}$ again as in \cite[Sec.~3.2.3]{KaKrVo20}.

Last, we fix the parameter $\gamma_B:=\frac{\delta}{3\,d}$, i.e., the failure probability of step~3 is bounded from above by $\gamma_B$, cf.\ \cite[Thm.~4.1]{Kae17}. Here, no modification is needed.

\subsubsection*{Final Step: Parameter insertion and union bounds}

The new parameter choices for $r$, $\gamma_{\operatorname{A}}$ and $\gamma_{\operatorname{B}}$ lead to the same failure probabilities for the detection of projected coefficients ($\frac{\delta}{3\,d\,s}$), Algorithm A for fixed $t$ and $i$ ($\frac{\delta}{3\,d\,s}$) and step~3 of Algorithm \ref{algo:sfft_general} ($\frac{\delta}{3\,d}$). Therefore, we can now use a union bound over the different steps of our algorithm similar to \cite[Thm.~9]{KaPoVo17}. It shows, that the total failure probability is now really bounded by terms less than $\delta$.

Finally, the sample complexity and computational complexity stated in Theorem \ref{thm:complexities} now follow directly using the above discussed choices $r = \lceil 2\,s\,\log(\frac{3\,d\,s\,G}{\delta})\rceil$, $\gamma_A:=\frac{\delta}{3\,d\,s\,G}$, and $\gamma_B:=\frac{\delta}{3\,d}$. The precise complexities for each step are given in Table \ref{table:proof2}.
\end{proof}

\begin{table}
\centering
\begin{tabular}{l l l}
& sample complexity & computational complexity \\
\hline
Step 1 & $d\,s\,N_\Gamma\,\log \frac{d\,s\,G}{\delta}$ & $d\,s\,G\,N_\Gamma\,\log^2 \frac{d\,s\,G\,N_\Gamma}{\delta}$ \\
\rule{0em}{1.25em}Step 2 (w.h.p.) & $d\,s\,\max(s,N_\Gamma)\,\log^2 \frac{d\,s\,G\,N_\Gamma}{\delta}$ & $d^2 \,s^2 \,G^2 \,N_\Gamma \,\log^3 \frac{d\,s\,G\,N_\Gamma}{\delta}$ \\
\rule{0em}{1.25em}Step 2 (w.c.) & $d\,s\,\max(s,N_\Gamma)\,\log^2 \frac{d\,s\,G\,N_\Gamma}{\delta}$ & $d^2 \,s^3 \,G^2 \,N_\Gamma \,\log^3 \frac{d\,s\,G\,N_\Gamma}{\delta}$ \\
\rule{0em}{1.25em}Step 3 & $\max(s\,G,N_\Gamma)\,\log \frac{d\,s\, G}{\delta}$ & $G\,\max(s\,G,N_\Gamma)\,\log \frac{d\,s\, G}{\delta}\, (d+\log(s\,G \,N_\Gamma))$
\end{tabular}
\caption{Same as Table \ref{table:proof1} but with the specifically chosen values for $r, \gamma_{\operatorname{A}}$ and $\gamma_{\operatorname{B}}$.}\label{table:proof2}
\end{table}

\begin{remark}\label{rem:thm}
The sample complexity of step~2 was the dominating term for the sFFT. Hence, we could neglect the sample complexity of step~3 there completely. In the usFFT, this sample complexity now contains a linear factor $G$, such that it is not neglectable for arbitrarily chosen $G$. However, if we can bound $G$ for example by $G \lesssim d\,s$, the sample complexity of step~3 is again asymptotically smaller than for step~2. Even more, since $G$ appears only in logarithmic terms of the sample complexity of step~2, we see, that the overall sample complexity of the usFFT is the same as for the sFFT in this case, i.e., the number of sampling locations is bounded in $\OO{d\, s\, \max(s, N_{\Gamma})\, \log^2 \frac{d\, s\, N_{\Gamma}}{\delta}}$ when assuming $G \lesssim d\,s$, cf.\ also \cite[Thm.~1.3]{KaKrVo20} for comparison. This is an important observation, since the amount of sampling locations is the crucial factor for the overall computational complexity of our algorithm due to the high computational cost of the underlying sampling algorithm, i.e., the PDE solver, as mentioned several times before.
\end{remark}

\section{Proof of Lemma \ref{lem_delta}}\label{app:proof_lem}

\begin{proof}
Since $i$ and $z_j$ in formula \eqref{eq:lattice_nodes} are integers, we know that $\tilde{y}_{i,j} \in \lbrace \frac{n}{M},\, n=0,...,M-1 \rbrace$ for all $i=0,...,M-1$ and $j=1,...,d$. In particular, since $M$ is prime and $z_{j_0} \not\equiv 0\imod{M}$, we have that $\lbrace\tilde{y}_{i,j_0}, i=0,...,M-1\rbrace = \lbrace \frac{n}{M},\, n=0,...,M-1 \rbrace$, so each $\frac{n}{M}$ is really attained at least once for some $i$ and $j$. Using this and the fact, that we are only considering $0 = \frac{0}{M} < \Delta < \frac{1}{2}\frac{1}{M} = \frac{1}{2M}$, we have
\begin{align*}
\min_{\substack{i=0,...,M-1\\j=1,...,d}} \left\lvert \tilde{y}_{i,j} - \Delta  \right\rvert &= \min_{n=0,...,M-1} \left\lvert \frac{n}{M} - \Delta  \right\rvert = \left\lvert\frac{0}{M} - \Delta\right\rvert = \Delta.
\end{align*}
On the other hand, we have
\begin{align*}
\min_{\substack{i=0,...,M-1\\j=1,...,d}} \left\lvert \tilde{y}_{i,j} - \left( \Delta + \frac12 \right) \right\rvert &= \min_{n=0,...,M-1} \left\lvert \frac{n}{M} - \left( \Delta + \frac12 \right) \right\rvert,
\end{align*}
where the minimum is attained for $n$ being the closest integer number to $M\Delta + \frac{M}{2}$. Since $0 < M\Delta < \frac{M}{2M} = \frac{1}{2}$ and $M$ odd, we conclude
\begin{align*}
\min_{n=0,...,M-1} \left\lvert \frac{n}{M} - \left( \Delta + \frac12 \right) \right\rvert &= \left\lvert \frac{M+1}{2M} - \left( \Delta+\frac{1}{2} \right) \right\rvert = \frac{1}{2M} - \Delta.
\end{align*}
Since the sum of these two minima is constant $\frac{1}{2M}$, we have the upper bound
\begin{align*}
\min \left\lbrace \min_{\substack{i=0,...,M-1\\j=1,...,d}} \left\lvert \tilde{y}_{i,j} - \Delta  \right\rvert, \min_{\substack{i=0,...,M-1\\j=1,...,d}} \left\lvert \tilde{y}_{i,j} - \left( \Delta + \frac12 \right) \right\rvert \right\rbrace = \min \left\lbrace \Delta, \frac{1}{2M} - \Delta \right\rbrace \leq \frac{1}{4M}.
\end{align*}
Finally, this upper bound is reached if and only if $\Delta = \frac{1}{4M}$ and hence
\begin{align*}
\Delta_{\text{opt}} = \argmax_{0 < \Delta < \frac{1}{2M}} \left\lbrace \min \left\lbrace \Delta, \frac{1}{2M} - \Delta \right\rbrace \right\rbrace = \frac{1}{4M}.
\end{align*}
\end{proof}

\begin{remark}
Note that the $\argmax$ in Lemma \ref{lem_delta} is not unique in general, since there also exist several values for $\Delta \geq \frac{1}{2M}$ attaining this maximum, e.g., $\Delta = \frac{3}{4M}$, which can be proven analogously. In our numerical experiments in Section \ref{sec:numerics}, we will always work with $\Delta_{\text{opt}} = \frac{1}{4M}$, which is the smallest optimal $\Delta > 0$ as we saw in the Theorem above.

Also, if we would neglect the assumption that $M$ is prime, we could run into problems if $z_j$ and $M$ are not coprime for all $j = 1,...,d$, since then $\lbrace\tilde{y}_{i,j}, i=0,...,M-1\rbrace$ is only a proper subset of $\lbrace \frac{n}{M},\, n=0,...,M-1 \rbrace$. But this case is neglectable, since our algorithm only uses prime lattice sizes $M$.
\end{remark}

\end{appendices}

{\scriptsize
\bibliographystyle{abbrv}

}

\end{document}